\documentclass[10pt, a4paper, oneside, reqno]{article}   	

\usepackage{geometry}                		
\usepackage[T1]{fontenc}
\usepackage[parfill]{parskip}    		
\usepackage{graphicx}				
\usepackage{fancyhdr}

\usepackage{amssymb}
\usepackage{mathtools}
\usepackage{amsmath}
\usepackage{dsfont}

\usepackage{newlfont}
\usepackage{amsfonts}
\usepackage{braket} 
\usepackage{subfig}
\usepackage{caption}
\usepackage{pdfpages}
\usepackage{enumerate} 
\usepackage{tikz} 
\usetikzlibrary{arrows.meta}
\usepackage{emptypage}

\usepackage{acronym}
\usepackage{textcomp}
\usepackage{siunitx}
\usepackage{verbatim}
\usepackage{eso-pic}

\usepackage{stackengine,graphicx}

\usepackage{indentfirst}

\usepackage{cite}

\usepackage{xcolor}
\pagecolor{white}

\usepackage[hidelinks]{hyperref} 

\usepackage{amsthm}

\newtheorem{thm}{Theorem}[section]

\newtheorem{prop}[thm]{Proposition}
\newtheorem{cor}[thm]{Corollary}
\newtheorem{lemma}[thm]{Lemma}

\theoremstyle{definition}

\numberwithin{equation}{section}

\stackMath

\newcommand{\abs}[1]{\left\lvert#1\right\rvert}
\DeclareMathSymbol\unlhd{\mathrel}{lasy}{"02}

\DeclareMathOperator{\sym}{sym}

\DeclareMathOperator{\supp}{supp}
\DeclareMathOperator{\rk}{rk}
\DeclareMathOperator{\id}{Id}
\DeclareMathOperator{\Hess}{H}
\DeclareMathOperator{\diam}{diam}
\DeclareMathOperator{\dist}{dist}

\let\phi\varphi
\let\epsilon\varepsilon
\newcommand{\ssubset}{\subset\joinrel\subset}

\newcommand*{\nsection}[1]{
    \section*{#1}
    \addcontentsline{toc}{section}{#1}
}

\usepackage{latexsym}

\usepackage{physics}

\usepackage{enumitem}

\usepackage{atbegshi,picture,afterpage}

\usepackage[affil-it]{authblk}
\usepackage[font=footnotesize,labelfont=bf]{caption}
\usepackage{floatrow}

\title{\vspace{-2cm}\textbf{Energy scaling for von K\'arm\'an elastic plates with positional constraints}}
\author{Giacomo Vizzari\thanks{Corresponding author. Email: \texttt{giacomo-vizzari@etud.univ-tln.fr}}}
\affil{Univ Toulon\\ Laboratoire de Mathématiques et Mathématiques appliquées IMATH\\ Toulon, France}
\date{October 28th, 2026}

\begin{document}

\pagenumbering{Roman} 
\maketitle

\begin{abstract}
    The von Kármán (vK) energy functional $\int \abs{\frac{1}{2}\nabla v\otimes \nabla v + \sym \nabla w}^2 + h^2 \abs{\nabla^2 v}^2 dx$, dependent on the bending stiffness $h$ of the plate, its vertical displacement $v$ and the horizontal one $w$, is non-convex. We make use of a relatively recent convex integration result from Lewicka and Pakzad and geometrical approximations (similar to the ones using origami maps) to find an upper bound for the scaling of the elastic energy as $h\to 0$. Through the interplay of the stretching and bending components of the energy we can attain a scaling of order $h^\frac{5}{3}$ for the upper-bound of a positional constraint problem. Furthermore, in doing so we will provide some new $\mathcal{C}^0$ approximation results for vK short maps through piecewise affine vK solutions.
\end{abstract}

\renewcommand\contentsname{}

\tableofcontents

\newpage

\pagenumbering{arabic}

\section{Introduction}
The ideas of this paper have to be considered in the broader topic that is the study of variational and geometrical problems related to the Föppl - von Kármán model for prestrained elastic plates (see \cite{lewicka2023calculus}). Specifically, we wish to study a system made of a thin film $\Omega\times[-\tilde{h}, \tilde{h}]$, which receives prestrain of some kind, i.e. an internal growth or shift. As such, the points composing our plate undergo some displacement $u$, which corresponds to the minimizer (if it exists) of some elastic energy functional. Specifically, the von Kármán (vK) ansatz states that to properly model the behavior of the elastic plate it is sufficient to study the mid-plane of the film (which we can identify with the 2-dimensional set $\Omega$), allowing us to reduce the dimensionality of the problem. G. Friesecke, R. D. James and S. Müller in \cite{FJM} have rigorously derived several non linear two dimensional von Kármán-like models through $\Gamma$-convergence theories, while studying the elastic energy scaling $E_h /h \sim h^\alpha, \alpha\ge 2$ in plates with several different boundary conditions or external forces. In case of a small strain, with the plate's crumpling differing only slightly from its resting position, we may look for some partial linearization of the problem, where the stretching energy functional is minimized for ${\epsilon_0}$-isometries of order 3 for some ${\epsilon_0}>0$. In fact, a different scaling can be considered for the in-plane component $w: \Omega\to \mathbb{R}^2$ with respect to the out-of-plane component $v:\Omega\to\mathbb{R}$ of the displacement as it follows:
\[
    u_{\epsilon_0}(x)=\id_2 (x)+{\epsilon_0} v(x)e_3+\epsilon_0^2w(x)
\]
where ${\epsilon_0}$ is the scale of the deformation. Compared with the usual models studied for example in \cite{contimaggi}, this deformation is a infinitesimal isometric immersion (more specifically an ${\epsilon_0}$-isometry of order 2) as long as 
\begin{equation}\label{vonkarmanintro}
    \frac{1}{2}\nabla v\otimes \nabla v + \sym \nabla w=0
\end{equation}
The absence of a linear term ${\epsilon_0} w(x)$ for the in-plane displacement leads to even better results: from straightforward computations we get that $\norm{\nabla u_{\epsilon_0}^T \nabla u_{\epsilon_0}-\id_2} = O(\epsilon_0^4)$ and as such the condition \eqref{vonkarmanintro} is sufficient to have an ${\epsilon_0}$-isometry of order 3. 


We want to consider a scenario of a plate positionally constrained within two layers and study its energy through the rescaled von Kármán elastic energy functional:
\begin{equation}\label{eq1}
    I^h[(w,v),\Omega]\coloneqq \int_\Omega \abs{\frac{1}{2}\nabla v\otimes \nabla v + \sym \nabla w}^2 + h^2\int_\Omega \abs{\nabla^2 v}^2
\end{equation}
where the first term corresponds to the \emph{stretching energy} and the second one to the \emph{bending energy} and the term $h^2$ here indicates some intrinsic resistance to bending deformation possessed by the plate, i.e. its (rescaled) thickness or bending stiffness. 

The model described by \eqref{eq1} can be ideally traced back to the Föppl--von Kármán energy of a plate of thickness $h$ normalized by its thickness: the membrane term scales like $h$ and the bending term like $h^3$ (see e.g.\cite{art3}, \cite{art2}, \cite{contimaggi}), and dividing by $h$ yields \eqref{eq1}, up to constants. The same functional can also be obtained, formally, from the von Kármán limit of \cite{multilayered} for a multilayered plate whose stiffness is concentrated in the inner plate of relative thickness $h$ (which is different from the thickness of De Benito Delgado and Schmidt, see figure 1): taking $W_0(t,F)=a_h(t)W(F)$ with $a_h(t)=h^{-1}a_0(t/h)$, $a_0$ even and compactly supported, one has $\int a_h=\int a_0$ and $\int a_h t^2=h^2\int a_0 s^2\,ds$, while the evenness of $a_h$ cancels the coupling term between stretching and bending. The choice $\tilde Q_2=|\cdot|^2$ is made for simplicity: any positive definite quadratic forms lead to the same scaling up to a constant. Finally, we add a positional constraint to the problem: we want the modelled plate to be confined between two rigid obstacles $v_{\inf}<v<v_{\sup}$, in the spirit of the confinement problem of \cite{contimaggi}.

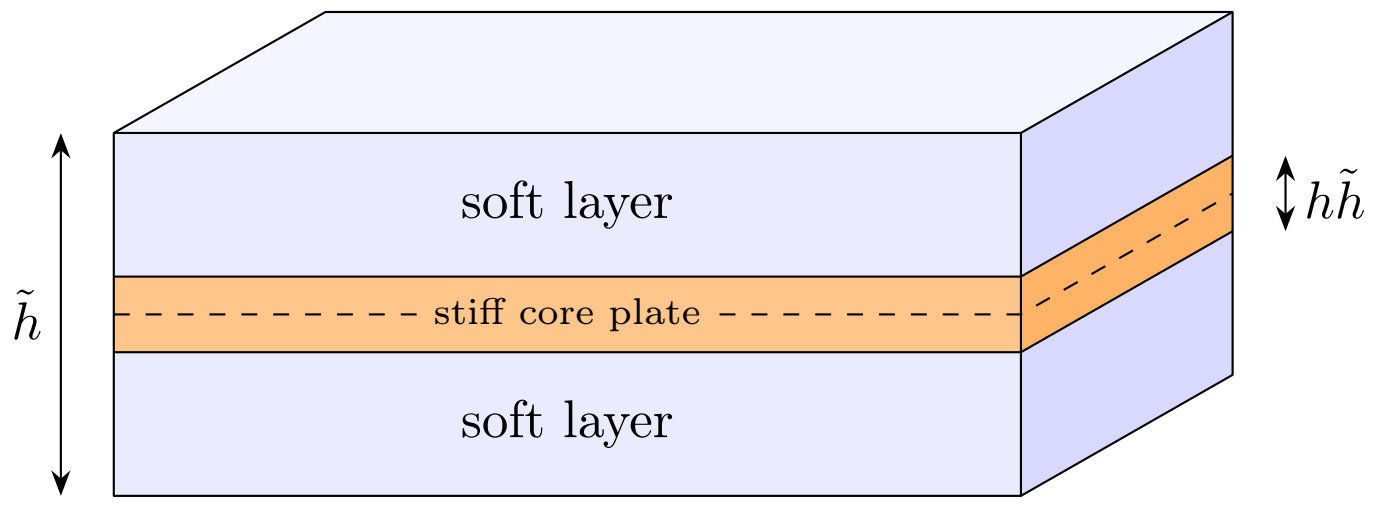
\begin{figure}[ht]\label{figure1}
\centering
\begin{tikzpicture}[>=Stealth, line join=round]
  \def\L{6}\def\H{2.4}\def\c{0.5}\def\dx{1.4}\def\dy{0.8}
  \pgfmathsetmacro{\lo}{\H/2-\c/2}
  \pgfmathsetmacro{\hi}{\H/2+\c/2}
  \pgfmathsetmacro{\mi}{\H/2}

  \fill[blue!8]    (0,0) rectangle (\L,\H);                              
  \fill[orange!45] (0,\lo) rectangle (\L,\hi);                           
  \fill[blue!15]   (\L,0) -- ++(\dx,\dy) -- ++(0,\H) -- (\L,\H) -- cycle; 
  \fill[orange!60] (\L,\lo) -- ++(\dx,\dy) -- ++(0,\c) -- (\L,\hi) -- cycle; 
  \fill[blue!4]    (0,\H) -- ++(\dx,\dy) -- ++(\L,0) -- (\L,\H) -- cycle; 

  \draw (0,0) rectangle (\L,\H);
  \draw (\L,0) -- ++(\dx,\dy) -- ++(0,\H) -- ++(-\L,0) -- (0,\H);
  \draw (\L,\H) -- ++(\dx,\dy);
  \draw (0,\lo) -- (\L,\lo) -- ++(\dx,\dy);
  \draw (0,\hi) -- (\L,\hi) -- ++(\dx,\dy);
  \draw[dashed] (0,\mi) -- (\L,\mi) -- ++(\dx,\dy);   

  \node at (\L/2,{(\hi+\H)/2}) {soft layer};
  \node at (\L/2,{\lo/2}) {soft layer};
  \node[fill=orange!45, inner sep=1pt] at (\L/2,\mi) {\scriptsize stiff core plate};

  \draw[<->] (-0.35,0) -- (-0.35,\H) node[midway,left] {$\tilde h$};
  \draw[<->] ({\L+\dx+0.35},{\lo+\dy}) -- ({\L+\dx+0.35},{\hi+\dy})
        node[midway,right] {$h\tilde h$};

\end{tikzpicture}
\caption{Sandwich plate of total thickness $\tilde h$ with a stiff core of
thickness $h\tilde h$ between two soft layers; the dashed line is the
midplane $\Omega$.}
\label{fig:sandwich}
\end{figure}

A variety of still open problems are related to this model: in fact, despite the simplifications, the energy functional is non-convex even when reduced to the midplane, which makes otherwise simple considerations not always attainable through standard methods. Furthermore, the different scaling of the stretching and bending energy makes the competition between these two terms non negligible. As such geometrical techniques are needed in order to derive, for example, any sort of bound on the scaling of the energy functional as the bending stiffness $h$ goes to 0.

Energy scaling problems for vK energy functionals have been tackled by several different studies, taking into account various kinds of geometrical and physical constraints, as for example \cite{art3} and \cite{art2}, studying the behavior of crumpling film under boundary conditions. The core concept is that the vK bending energy in \eqref{eq1} naturally scales like $h^2$, while the stretching energy has a lower order scaling of $O(1)$ in $h$. Hence to get the best possible scaling for both, we would like to minimize the stretching energy as much as possible, which, as we have seen in \eqref{vonkarmanintro}, is attained for a linearized isometric immersion, and sacrifice the bending term in turn. The idea comes from the study on the nonlinear FvK model with positional constraints by S. Conti and F. Maggi in \cite{contimaggi}: they considered the scaling of the elastic energy functional for a thin film, as it was crumpled inside a ball of small radius $r$, finding an upper bound of $h^\frac{5}{3}$ (and some partial results for verifying it as a lower bound as well). To do that, they have utilized a fully geometrical approach by "origami maps": first, by the Nash-Kuiper theorem on isometric immersions they approximate the $W^{1,\infty}$ \emph{short maps} with piecewise affine isometries obtained by folding the sheet along some edges (origami maps); then they  proceed to approximate the piecewise affine map with a smooth one away from the vertices of their triangulations, following the ridge scaling ideas of A. E. Lobkovsky and T. A. Witten in \cite{physics} (one could also look at the results obtained by P. Gladbach and H. Olbermann with a convexity assumption for FvK cones and ridges in \cite{cones2} and \cite{cones1}; or at the full bound for the ridge construction depending on thickness and bending angle for the FvK setting obtained by S. Venkataramani in \cite{venkataramani}).

The constant parallel between the vK and FvK energy functional leads to believe that a similar problem with positional constraints can attain the same scaling of $h^\frac{5}{3}$ with fairly similar methods, but with some major differences. On the one hand, our nontrivial in-plane displacement $w$ for the plate will lead to different reparametrizations and approximations. On the other hand, dealing with non-exact isometries will mean that we will not be able to use euclidean distance to formulate geometrical conditions; hence we had to find equivalent ideas (see \ref{preplemma2}) and adapt them to the desired construction. Finally, a parallel to the Nash-Kuiper theorem for isometric embeddings from \cite{Nash-Kuiper} used by Conti and Maggi, a theorem proved by M. Lewicka and M. R. Pakzad (see \cite{mongeampere4}) provided an analogous result for linearized isometric immersions, which are strictly related to our energy functional. By using that, the following upper bound will be proven in this paper.

\begin{thm}\label{mainthm}
Given $\Omega\subset \mathbb{R}^2$ bounded Lipschitz domain, $(w_0, v_0)\in \mathcal{C}^1(\bar{\Omega};\mathbb{R}^2\times\mathbb{R})$ such that
\begin{equation}\label{short}
    \frac{1}{2} \nabla v_0 \otimes \nabla v_0 + \sym \nabla w_0 \le 0
\end{equation}
then for all $r > 0$ there exist $h^*=h^*(r, (w_0, v_0))>0$ and $C^*=C^*(\Omega, r, (w_0, v_0))$ such that, for all $h<h^*$, there exists some $(w_h, v_h)\in \mathcal{C}^2(\Omega, \mathbb{R}^2\times\mathbb{R})$ such that $\norm{(w_h, v_h)-(w_0,v_0)}_{L^{\infty}(\Omega; \mathbb{R}^2\times\mathbb{R})}<r $ and
\[
    I^h[(w_h, v_h), \Omega]\le C^* h^{\frac{5}{3}}
\]
\end{thm}

We can see that this scaling is not intuitive to attain. In fact, Theorem \ref{LPthm} provides exact solutions $(w,v)\in C^{1,\beta}$ uniformly close to $(w_0,v_0)$. Mollifying them at scale $\sigma$, we would ultimately get $I^h\lesssim \sigma^{4\beta}+h^2\sigma^{2\beta-2}$, and the choice $\sigma=h^{1/(1+\beta)}$ gives $I^h\lesssim h^{4\beta/(1+\beta)}$ (cf.\ the analogous bound for prestrained films in \cite{thm1.11}, Theorem 1.11). For the best known exponent $\beta<1/3$ (see \cite{frac13}) this is worse than $h$, and reaching $h^{5/3}$ would require $\beta\ge5/7>2/3$ which is beyond the expected rigidity threshold. The scaling of Theorem \ref{mainthm} is therefore not accessible by regularizing convex integration solutions. On the other hand, if we were to have some $(w_j,v_j)$ such that $I^{h_j}[(w_j, v_j), \Omega]\le C^* h_{j}^{2}$ for all $j$ and a fixed $(w_0,v_0)$ with $\norm{(w_j, v_j)-(w_0,v_0)}_{L^{\infty}}\le r$, then it can be proven that up to a subsequence $v_{h_j}$ weakly converges in $W^{2,2}$ to some $v$ such that $\det \nabla ^2 v=0$, and thus developable on convex domains by \cite{developability} (result that can be also extended to non-convex domains, see \cite{developability2} Cor 2.2); this would necessarily lead to a contradiction along the segments on which $v$ is affine for small values of $r$ when $v_0$ is for example quadratic (e.g. $v_0(x)=\abs{x}^2$ and $w_0(x)=-2x$ on the unit ball we necessarily have $\norm{v_0-v}_{L^\infty}\ge \frac{1}{8}$). Hence, necessarily for $r<1/8$: 
\[
	\lim_{h\to 0}\frac{1}{h^2}\inf\{I^h[(w,v),\Omega]\;:\; |(w,v)-(w_0,v_0)|_{L^\infty}< r\}=+\infty
\]
Meanwhile our result for $h^{5/3}$ would be consistent with this newfound threshold, however proving its optimality is left open.

We may see how the condition \eqref{short} gives us a infinitesimal order "\emph{short map}" to approximate to prove this theorem; as such, it will be helpful to work with the equivalent of \emph{origami maps} for our setting as well. For the definition of triangulations and origami maps we will refer to \cite{contimaggi}, and for convenience of the reader report that same definitions here, modified for our problem. Given a bounded open $\Omega\subset \mathbb{R}^2$ we will define a triangulation $\mathcal{T}=\{T_j\}_j$ as a finite family of closed, nondegenerate triangles with disjoint interiors such that $\Omega$ is contained in $\Omega_\mathcal{T}\coloneqq \bigcup_j T_j$ and such that each pair $(T_i, T_j)$ is either disjoint, shares a corner or shares a side. On such a triangulation we will define an \emph{origami map for the vK constraint} as $(\mathcal{T}, (w,v))$, with $(w,v)\in W^{1,\infty}(\Omega_\mathcal{T};\mathbb{R}^3)$ affine on each $T_i\in \mathcal{T}$, and such that it satisfies \eqref{vonkarmanintro}. To every triangulation we associate its degeneracy factor (which we usually want to bound in a proper geometrical construction):
\[
    \Gamma(\mathcal{T})\coloneqq \max\left\{l_{\max{}},\nu,\frac{1}{\tau_{\min{}}}\right\}
\]
where $l_{\max{}}$ is the maximum side length of the triangles in $\mathcal{T}$, $\tau_{\min{}}$ is their their minimum slope and $\nu$ is the maximum number of triangles of $\mathcal{T}$ that meet at a point.

The proof of Theorem \ref{mainthm} can be divided in two main parts. Firstly, by using a Lewicka and Pakzad theorem from \cite{mongeampere4} we will restrict ourselves to $\mathcal{C}^{1,\beta}$ class functions solving \eqref{vonkarmanintro}. Then via triangulation techniques we will provide an approximation of such a function with an \emph{origami map}. This first part is covered in Section 2. Specifically, Section 2 is the core of our construction: in it we attain novel and interesting results for an origami construction for solutions to vK constraint without a Euclidean distance. 

Secondly, we will approximate the origami map to be able to compute the bending energy over it and provide a bound. This is necessary since origami maps are piecewise affine and the bending component of the energy functional \eqref{eq1} "explodes" over the edges of the triangulation. The idea is to approximate the map around the sides of our triangles in regions that "shrink" towards the vertices. The area around the vertices will be cut off to avoid singularities and will be treated separately. This construction is covered in Section 3.

Finally, we may see that a similar obstacle problem result to Theorem 1.4 from \cite{contimaggi} can be consequently obtained.

\begin{cor}\label{positional}
    Let $\Omega\subset \mathbb{R}^2$ be a bounded Lipschitz domain, $v_{\sup}$ and $v_{\inf}$ be two functions in $\mathcal{C}^0(\bar{\Omega};\mathbb{R})$ such that $v_{\sup}-v_{\inf}>r$ for some $r>0$. Let $h^*=h^*(r, v_{\sup}, v_{\inf})>0$, then there exist $C^*=C^*(\Omega, r, v_{\sup}, v_{\inf})$ such that for all $0<h<h^*$ we can find $(w_h, v_h)\in \mathcal{C}^2(\Omega, \mathbb{R}^2\times\mathbb{R})$ that satisfies
    \[
    I^h[(w_h, v_h), \Omega]\le C^* h^{\frac{5}{3}}
    \]
    and $v_{\inf}<v_h<v_{\sup}$ on $\Omega$.
\end{cor}

In this positioning problem we are intuitively restricting our out-of-plane deformation $v_h$ to be fully contained in a section of $\Omega\times \mathbb{R}$ delimited by two obstacles $v_{\inf}$ and $v_{\sup}$. It is intuitively a similar kind of result to the "paper crumpled inside a ball" idea of the above mentioned Theorem from Conti and Maggi. However, it is important to notice how this is not attainable by direct application of Conti and Maggi to $u_{\epsilon_0}$ since their estimates have no control over the $\epsilon_0^2$-order terms of $\norm{\nabla u_{\epsilon_0}^T \nabla u_{\epsilon_0}-\id_2}$ and furthermore the constants they have found would depend from $\epsilon_0$. Hence, Theorem \ref{mainthm} and its Corollary cannot be directly deduced from \cite{contimaggi}.

%
%
%
%
%

\section{Construction of an affine solution to the von Kármán constraint} 

We begin by focusing on the construction of an origami map for the $vK$ constraint that approximates our starting function. To do so we start from a \textit{subsolution} to the von Kármán constraint \eqref{vonkarmanintro}, i.e. a function $(w_0,v_0)\in \mathcal{C}^1(\bar{\Omega};\mathbb{R}^2\times \mathbb{R})$ that satisfies the constraint:
\begin{equation}\label{subsolution}
    \frac{1}{2} \nabla v_0 \otimes \nabla v_0+\sym \nabla w_0 \le 0
\end{equation}
We aim to obtain the following:
\begin{thm}\label{mainthm3}
    Every short vK map can be approximated in $\mathcal{C}^0$ by piecewise affine solutions on a non-degenerate triangulation.
    
    In other words, let $(w_0,v_0)\in \mathcal{C}^1(\bar{\Omega};\mathbb{R}^2\times \mathbb{R})$, on some $\Omega\subset \mathbb{R}^2$ bounded Lipschitz domain, such that it satisfies \eqref{subsolution}. Then for every $\epsilon >0$ there exists some piecewise affine map $(\tilde{w}_\epsilon, \tilde{v}_\epsilon)\in W^{1,\infty}(\Omega;\mathbb{R}^2\times\mathbb{R})$ satisfying \eqref{vonkarmanintro} on a triangulation $\mathcal{T}_\epsilon$ of $\Omega$ such that
    \begin{align*}
        &\lVert \tilde{v}_\epsilon-v_0\rVert_{L^\infty(\Omega ; \mathbb{R})}\le \epsilon \quad
        &\lVert \tilde{w}_\epsilon-w_0\rVert_{L^\infty(\Omega ; \mathbb{R}^2)}\le \epsilon \\
        &\norm{\nabla \tilde{v}_\epsilon}_{L^\infty} \le C_\epsilon(\norm{(w_0,v_0)}_{\mathcal{C}^1}) \quad
        &\lambda_\epsilon \le C_\epsilon(\norm{(w_0,v_0)}_{\mathcal{C}^1})
    \end{align*}
    where $\lambda_\epsilon$ is the norm of the maximum jump of $\nabla \tilde{v}_\epsilon$ in $\Omega$. Furthermore we may have $\Gamma(\mathcal{T}_\epsilon)\le C$ for some universal constant $C$.
\end{thm}

\subsection{Structure of the proof}
Before providing the proof we will start by looking at some preparatory lemmas. We will take the first one directly from \cite{contimaggi} (Lemma 4.1).

\begin{lemma}\label{preplemma1}
    Let $\Omega \subset \mathbb{R}^n$ be a bounded Lipschitz domain. Then there is a sequence of diffeomorphisms $\Phi_k\in \mathcal{C}^\infty (\mathbb{R}^n,\mathbb{R}^n)$ such that
    \[
        \lim_{k\to\infty}\sup_{p\in \mathbb{R}^n} [\abs{\Phi_k(p)-p}+\abs{\nabla \Phi_k(p)-\id_n}+\abs{\Hess_{\Phi^{-1}_k(p)}}]=0
    \]
    and $\Omega\ssubset \Phi_k (\Omega)$.
\end{lemma}
The proof can be found in \cite{contimaggi}, with the only added requirement here being the uniform convergence to 0 of the norm of the Hessian matrix of $\Phi^{-1}_k$, which can be directly deduced from the proposed construction for $\Phi_k$ and the inverse function theorem.

\begin{lemma}\label{preplemma2}
Let $\Omega\subset \mathbb{R}^2$ be an open bounded domain.

(i) Let $(w,v)\colon \Omega \to \mathbb{R}^2\times\mathbb{R}$ be a mapping at least piecewise of class $\mathcal{C}^{1,\beta}$ and such that 
\begin{equation}\label{deltaeq}
    \frac{1}{2} \nabla v \otimes \nabla v+\sym \nabla w = -\delta\id_2 
\end{equation}
for some $\delta\in \mathbb{R}$, then we may prove that 
for any pair of points $x,y\in \Omega$ such that $\gamma(t)\coloneqq x+t(y-x)\subset \Omega$ for all $t\in[0,1]$ we have:
\begin{equation}\label{delta}
    \frac{1}{2} \abs{v(x)-v(y)}^2+(w(x)-w(y))\cdot(x-y) \le -\delta \abs{x-y}^2
\end{equation}
and, more specifically, if $v$ is affine then we have an equality 
in \eqref{delta}.

(ii) Let $(w,v)$ be as in (i), $\kappa \in \left(1,\frac{3}{2}\right)$, then there exists some $l^*=l^*(\kappa,\delta, \norm{(w,v)}_{\mathcal{C}^{1,\beta}})>0$ small enough such that for any pair of points $x,y\in \Omega$ such that $\gamma(t)\coloneqq x+t(y-x)\subset \Omega$ for all $t\in[0,1]$, and $\abs{x-y}<l^*$ we have:
\begin{equation}\label{deltaleft}
    -\kappa\delta \abs{x-y}^2\le 
    \frac{1}{2} \abs{v(x)-v(y)}^2+(w(x)-w(y))\cdot(x-y)
\end{equation}

(iii) Conversely, if $(w,v):\Omega\to \mathbb{R}^2\times\mathbb{R}$ is an affine map such that for any pair of points $x,y\in \Omega$ \eqref{delta} holds, then we have:
\begin{equation}\label{delta2}
    \frac{1}{2} \nabla v \otimes \nabla v+\sym \nabla w \le -\delta\id_2 
\end{equation}
where if \eqref{delta} holds with an equality, then \eqref{delta2} is an equality as well.
\end{lemma}
This lemma will be extensively used in the next few sections, since it allows us to change our starting condition \eqref{subsolution} to one that can be simply verified on segments. Specifically if $\delta=0$, it suffices to have: 
\begin{equation}\label{cond}
\frac{1}{2} \abs{v(x)-v(y)}^2+(w(x)-w(y))\cdot(x-y) =0
\end{equation}
for any pair of points $x,y\in \Omega$ such that $\abs{x-y}<l$ to have that an affine (or piecewise affine) map $(w,v)$ solves \eqref{vonkarman}.

\begin{proof} 
(i) Let us assume that $(w,v)$ is of class $\mathcal{C}^{1,\beta}$ (the piecewise case is analogous and only requires us to split the integration along $\gamma$) and let $\xi=x-y \in \mathbb{R}^2$. By applying it and its transposed to both sides of \eqref{deltaeq} we get:
\[
    \frac{1}{2}(\nabla v\cdot \xi)^2+\xi\cdot\nabla w \xi = -\delta\abs{\xi}^2
\]
and by integrating along $\gamma$, by the Fundamental Theorem of Calculus:
\begin{equation}\label{deltaproof}
    \frac{1}{2}\int_0^1(\nabla v(\gamma(t))\cdot \xi)^2dt+\xi\cdot (w (x)-w(y)) = -\delta\abs{x-y}^2
\end{equation}
On the other hand we may see that:
\[
    \frac{1}{2} \abs{v(x)-v(y)}^2+(w(x)-w(y))\cdot(x-y)=\frac{1}{2}\left(\int_0^1\nabla v(\gamma(t))\cdot \xi dt\right)^2+\xi\cdot (w (x)-w(y)) 
\]
we may combine that with \eqref{deltaproof} by using Jensen's inequality to get \eqref{delta}. 

Moreover, if $v$ is affine then $\nabla v$ is constant and the equality condition for Jensen is satisfied. 

(ii) Let us expand by Taylor, we have:
\begin{gather*}
    v(x) - v(y) = \nabla v(y) \cdot (x-y) + o(|x-y|)\\
    w(x) - w(y) = \nabla w(y) (x-y) + o(|x-y|)
\end{gather*}
Hence:
\begin{gather*}
    (v(x) - v(y))^2 = |\nabla v(y) \cdot (x-y)|^2 + o(|x-y|^2)\\
    2(w(x) - w(y)) \cdot (x-y) = 2(x-y)^T \nabla w(y) (x-y) + o(|x-y|^2)
\end{gather*}

Summing them up gives us:
\[
(x-y)^T \left[ \nabla v \otimes \nabla v + \nabla w + \nabla w^T \right] (x-y) + o(|x-y|^2)=-2\delta |x-y|^2 + o(|x-y|^2)
\]
By definition of $o(|x-y|^2)$ for any $\epsilon_\delta > 0$, there exists an $l^*$ such that for $|x-y| < l^*$, the error is smaller than $\epsilon_\delta |x-y|^2$. If we choose $\epsilon_\delta = (\kappa-1) \delta$ for some $\kappa > 1$, then for $|x-y| < l^*=l^*(\kappa,\delta)$
\[
-\delta |x-y|^2 + o(|x-y|^2)\ge-\kappa  \delta |x-y|^2
\]
which proves our thesis.

(iii) If we define $v(x)=a\cdot x+c, w(x)=Bx+C$ for some $c\in\mathbb{R}, a,C\in \mathbb{R}^2, B\in \mathbb{R}^{2\times 2}$. Let $x\in \Omega$, then let us take $\xi\in\mathbb{R}^2$ as any vector such that $\norm{\xi}\le \dist(x, \partial\Omega)$; then we have that there exists $y\in\Omega$ such that $\xi=x-y$ and:
\[
\begin{split}
    \frac{1}{2} \abs{v(x)-v(y)}^2+(w(x)-w(y))\cdot(x-y)&=\frac{1}{2} (a\cdot \xi)^2+ \xi\cdot B \xi \\
    &=\xi^T\left(\frac{1}{2} a \otimes a+\sym B\right)\xi
\end{split}
\]
On the other hand by \eqref{delta}:
\[
	 \frac{1}{2} \abs{v(x)-v(y)}^2+(w(x)-w(y))\cdot(x-y)\le \xi^T(-\delta \id)\xi
\]
then \eqref{delta2} derives directly from the arbitrary choice of $\xi$. The equality case is analogous.
\end{proof}

Finally, we will look at the following lemma which guarantees that \eqref{cond} is consistent with affine interpolation:

\begin{lemma}\label{interpolation}
Let $(w,v)$ be an affine function on a triangle $[x,y,z]$ such that \eqref{cond} holds for any pair of vertices of the triangle. Then \eqref{cond} holds for any pair of points in $[x,y,z]$.

\end{lemma}

\begin{proof}
To prove that for any pair of points $p,q$ in $[x,y,z]$ the thesis holds, it suffices to prove that 
\begin{equation}\label{interpol}
	\Delta(x,p)\coloneqq\frac{1}{2} \abs{v(x)-v(p)}^2+(w(x)-w(p))\cdot(x-p)= 0
\end{equation}
holds for any pair consisting of a point $p$ inside the triangle and a vertex $x$. In fact, by connecting $p$ with the vertices of the triangle we may split it into three smaller triangles. Let $[p,x,y]$ be, without loss of generality, the triangle containing $q$. By \eqref{interpol} we have that the hypothesis of the Lemma hold on $[p,x,y]$, hence \eqref{interpol} holds for the vertex $p$ and the inner point $q$.

Let us first verify that if $\Delta(x,y)= 0$, then we have that $\Delta(x,p)=0$ for any point $p= \rho x + (1-\rho) y$ for $\rho\in (0,1)$ in the segment $[x,y]$:
\begin{align*}
        \Delta(x,p)=&\frac{1}{2} \abs{v(x)-v(\rho x+(1-\rho) y)}^2+(w(x)-w(\rho x+(1-\rho) y))\cdot(x-(\rho x+(1-\rho) y))\\
        =&\frac{1}{2} \abs{v(x)-\rho v(x)-(1-\rho)v( y)}^2+(w(x)-\rho w( x)-(1-\rho)w(y))\cdot(x-\rho x-(1-\rho) y))\\
        =&\frac{1}{2} (1-\rho)^2\abs{v(x)-v(y)}^2+(1-\rho)^2(w(x)-w(y))\cdot(x-y)=0
\end{align*}

Then, given a point $p=\rho y + (1-\rho) z$, for some $\rho\in (0,1)$, in the segment $[y,z]$ we want to prove that $\Delta(x,p)=0$. This is indeed sufficient to prove our thesis, since we have then proven that $\Delta(x,p')=0$ for any inner point $p'$ of the segment $[x,p]$ and thus, given the arbitrary choice of $p\in [y,z]$, the same holds for any inner point of the triangle.

\begin{align*}
        \Delta(x,p)=&\frac{1}{2} \abs{v(x)-v(\rho y+(1-\rho) z)}^2+(w(x)-w(\rho y+(1-\rho) z))\cdot(x-(\rho y+(1-\rho) z))\\
        =&\rho^2 \Delta(x,y) + (1-\rho)^2 \Delta(x,z) + \rho(1-\rho) [(v(x)-v(y))(v(x)-v(z))+\\&+(w(x)-w(y))\cdot(x-z)+(w(x)-w(z))\cdot(x-y)]
\end{align*}
Since $\Delta(x,y)=0=\Delta(x,z)$, it remains to be proven that 
\begin{equation}\label{partialinterp}
	(v(x)-v(y))(v(x)-v(z))+(w(x)-w(y))\cdot(x-z)+(w(x)-w(z))\cdot(x-y)=0
\end{equation}
But looking at the first term of this equation $V=(v(x)-v(y))(v(x)-v(z))$, and remembering that $\Delta(x,y)=\Delta(x,z) =\Delta(y,z)=0$, we have that:
\begin{align*}
	V&=\frac {1}{2} \left[-(v(z)-v(y))^2+(v(x)-v(y))^2+(v(x)-v(z))^2\right]\\
	&= (w(y)-w(z))\cdot(y-z)-(w(x)-w(y))\cdot(x-y)-(w(x)-w(z))\cdot(x-z)\\
	&=(w(y)-w(x))\cdot (y-z + x-y)+(w(x)-w(z))\cdot(y-z+z-x)\\
	&=-(w(x)-w(y))\cdot (x-z)-(w(x)-w(z))\cdot(x-y)
\end{align*}
and by substituting into \eqref{partialinterp} we have our thesis.
\end{proof}

Let us now turn to the proof of Theorem \ref{mainthm3}. The idea (which is fundamentally geometrical) is the following: we would like to modify pointwise our displacement function on some triangulation such that we can obtain an origami-like map that satisfies the vK constraint via an affine interpolation. To do so, we will need to first enlarge our domain $\Omega$ via a reparametrization to avoid problems near $\partial \Omega$. 

To do so, we will first restrict ourselves to class $\mathcal{C}^{1,\beta}$ solutions of \eqref{vonkarmanintro} via a result from Lewicka and Pakzad (see \cite{mongeampere4}, Theorem 1.2). This theorem is an equivalent of the classical Nash-Kuiper's theorem for $\mathcal{C}^1$ isometric embeddings (see \cite{Nash-Kuiper}) applied instead to linearized isometric immersions:

\begin{thm}\label{LPthm}
    Let $\Omega\subset \mathbb{R}^2$ be an open bounded domain. Let $v_0\in \mathcal{C}^1 (\bar\Omega;\mathbb{R})$ and $w_0\in \mathcal{C}^1 (\bar\Omega;\mathbb{R}^2)$ and $A_0\in \mathcal{C}^{0,\alpha}(\bar{\Omega}, \mathbb{R}_{sym}^{2\times 2})$ for some $\alpha\in (0,1)$, be such that:
    \[
        \exists c_0>0\quad\quad A_0 -\left(\frac{1}{2}\nabla v_0 \otimes \nabla v_0 + \sym \nabla w_0\right)>c_0 \id_2\quad \quad \textrm{in }\bar{\Omega} 
    \]
    Then, for every exponent $\beta$ in the range:
    \[
        0<\beta <\min\left\{\frac{1}{7},\frac{\alpha}{2}\right\}
    \]
    there exist sequences $v_n\in \mathcal{C}^{1,\beta}(\bar{\Omega};\mathbb{R})$ and $w_n\in \mathcal{C}^{1,\beta}(\bar{\Omega};\mathbb{R}^2)$ which converge uniformly to $v_0$ and $w_0$ respectively, and which satisfy:
    \[
        A_0=\frac{1}{2}\nabla v_n \otimes \nabla v_n + \sym \nabla w_n \quad\quad \textrm{in }\bar\Omega
    \]
\end{thm}

The proof can be found in \cite{mongeampere4}. The result has been improved by W. Cao and L. Székelyhidi Jr in \cite{1/5first} and M. Lewicka in \cite{1/5lewi} up to $1/5$ and later on by W. Cao, J. Hirsch and D. Inauen who attained $\beta <\min\left\{\frac{1}{3},\frac{\alpha}{2}\right\}$ in \cite{frac13}. We may also remark that via (3.3) and (3.16) of \cite{mongeampere4} we may have:
\begin{equation}\label{gradient1}
	\norm{(\nabla w_n,\nabla v_n)}_{\mathcal{C}^0}\le C=C(\norm{w_0,v_0}_{\mathcal{C}^1})
\end{equation}

In the next step, we will use a rescaling to get a $\mathcal{C}^{1,\beta}$ map that satisfies \eqref{deltaeq} for some $\delta>0$. We will now have a definite gap in the vK constraint which will give us the ability to actually modify our function and finally adjust it into a piecewise affine solution to \eqref{vonkarmanintro}. The piecewise affine construction in section 2.2 is the focal part of the proof here, and will be made via affine extension after defining the function on a limited set of vertices. In the following proof we will omit that construction to see it more in detail in the following section.

\begin{proof}[Proof of Theorem \ref{mainthm3}]
    For all $k\in \mathbb{N}$, $k$ large enough, let $\Omega_k\coloneqq\Phi_k(\Omega)$ where $\Phi_k$ is a diffeomorphism given by Lemma \ref{preplemma1}, with $\Omega\ssubset\Omega_k$. Let us define a sequence $(w_k, v_k)\in \mathcal{C}^1(\bar{\Omega}_k;\mathbb{R}^2\times \mathbb{R})$ in the following way:
    \begin{gather*}
        v_k(x)= v_0(\Phi^{-1}_k (x))\\
        w_k(x)= \nabla\Phi^{-T}_k w_0(\Phi^{-1}_k (x))-\mu_k x\\
        \mu_k=\norm{ \Hess_{\Phi_{k,1}^{-1}} w_{0,1}+\Hess_{\Phi_{k,2}^{-1}} w_{0,2}}_\infty +\frac{1}{k}
    \end{gather*}

    Where $\Hess_\Phi$ is the hessian matrix of $\Phi$ and we have used the following notation to differentiate between the components of our functions: $(w_{0,1},w_{0,2})=w_0$ and $(\Phi_{k,1},\Phi_{k,2})=\Phi_k$ for all $k\ge 1$. By Lemma \ref{preplemma1} $\norm{v_k-v_0}_{L^\infty(\Omega;\mathbb{R})}\to 0$ and $\norm{w_k-w_0}_{L^\infty(\Omega;\mathbb{R}^2)}\to 0$ for $k\to \infty$. Furthermore, through some straightforward tensor calculations and with a proper change of coordinates we get:
    \[
    \begin{split}
        \frac{1}{2} \nabla v_k\otimes \nabla v_k + \sym \nabla w_k=&\nabla\Phi^{-T}_k \left[\left( \frac{1}{2} \nabla v_0\otimes \nabla v_0 + \sym \nabla w_0\right)\right]\nabla\Phi^{-1}_k+\\
        &+\Hess_{\Phi_{k,1}^{-1}}w_{0,1} +\Hess_{\Phi_{k,2}^{-1}} w_{0,2}-\norm{ \Hess_{\Phi_{k,1}^{-1}} w_{0,1}+\Hess_{\Phi_{k,2}^{-1}} w_{0,2}}_\infty \id_{2}-\frac{1}{k}\id_{2}
    \end{split}
    \]
    Hence $(w_k,v_k)$ satisfy \eqref{delta2} in the form:
    \[
     \frac{1}{2} \nabla v_k\otimes \nabla v_k + \sym \nabla w_k\le -\frac{1}{k} \id_2
     \]
     Now we may apply Theorem \ref{LPthm} for $A_0=0$ to show that for each $k$ there exists a sequence $(w_k^j, v_k^j)\in \mathcal{C}^{1,\beta}(\Omega_k; \mathbb{R}^2\times \mathbb{R})$, for some $\beta>0$, that converges uniformly to $(w_k,v_k)$ for $j\to \infty$ and such that $(w_k^j, v_k^j)$ solves \eqref{vonkarmanintro} for each choice of $k,j \ge 1$. By some kind of diagonal selection we can then get a sequence of class $\mathcal{C}^{1,\beta}$ solutions to \eqref{vonkarmanintro} that uniformly converge to $(w_0, v_0)$ on $\Omega$. For simplicity of notation we will from now on consider a $(w_k^k, v_k^k)$ in this diagonally selected sequence such that 
    \[
        \lVert (w_k^k,v_k^k)-(w_0,v_0)\rVert_{L^\infty(\Omega ; \mathbb{R}^2\times\mathbb{R})}\le \frac{\epsilon}{2},
    \]
    and rename $(w, v)\coloneqq(w_k^k, v_k^k)$, $\Omega^*\coloneqq \Omega_k$. Let us show that we can provide $(\tilde{w}_\epsilon,\tilde{v}_\epsilon)$ origami maps for the vK constraint (i.e. piecewise affine solutions to \eqref{vonkarmanintro}) with maximum jump $\lambda_\epsilon$ of $\nabla v_\epsilon$ and such that:
    \[
        \lVert (w,v)-(\tilde{w}_\epsilon,\tilde{v}_\epsilon)\rVert_{L^\infty(\Omega ; \mathbb{R}^2\times \mathbb{R})}\le \frac{\epsilon}{2}, \quad
        \norm{\nabla\tilde{v}_\epsilon}_{L^\infty}\le C_\epsilon(\norm{(w,v)}_{\mathcal{C}^1})
    \]
    which would prove Theorem \ref{mainthm3}.

    Let $\delta\in (0,1/2)$ to be fixed later on and let us now provide a rescaling of $(w, v)$ on $\Omega^*$ by a factor $1-\delta$. We define:
    \begin{gather*}
        v_\delta (x)\coloneqq (1-\delta) v(x)\\
        w_\delta(x)\coloneqq(1-\delta)^2(w(x)+x)-x
    \end{gather*}
    It follows that $(w_\delta, v_\delta)$ are still of class $\mathcal{C}^{1,\beta}$ for $\beta>0$ and:
    \[
        \frac{1}{2} \nabla v_\delta \otimes \nabla v_\delta + \sym \nabla w_\delta =(1-\delta)^2 \left(\frac{1}{2} \nabla v \otimes \nabla v + \sym \nabla w\right)+ (\delta^2-2\delta) \id_2=-(2-\delta)\delta \id_2
    \]
    We want to study $(w_\delta, v_\delta)$ at a scale where the map is almost affine to have proper bounds on the gradient jumps. So, given $\kappa\in(1,3/2)$ by Lemma \ref{preplemma2} we have that, for any two points $x,y\in \Omega^*$ such that $\gamma(t)\coloneqq x+t(y-x)\subset \Omega^*$ for all $t\in[0,1]$ and $\abs{x-y}<l^*$ for some $l^*=l^*(\kappa,\delta, \norm{(w_\delta,v_\delta)}_{\mathcal{C}^{1,\beta}})>0$, we get:
\begin{equation}\label{delta3}
    -\kappa(2-\delta)\delta \abs{x-y}^2\le 
    \frac{1}{2} \abs{v_\delta(x)-v_\delta(y)}^2+(w_\delta(x)-w_\delta(y))\cdot(x-y) \le -(2-\delta)\delta \abs{x-y}^2
\end{equation}
    Now let us choose:
    \[
        l\le \min\left\{l^*(\kappa,\delta)/2,\delta, \dist(\Omega,\partial\Omega^*)/2  \right\}
    \]
    We start by considering a partition of $\mathbb{R}^2$ in equilateral triangles of side length $l$ and we let $\mathcal{T}_l$ be the sub-triangulation defined by those which are contained in $\Omega^*$. By the choice of $l$ we have that $\Omega \subset \Omega_{\mathcal{T}_l}\subset \Omega^*$ (with $\Omega_{\mathcal{T}_l}$ domain formed by the union of all the triangles of $\mathcal{T}_l$), also \eqref{delta3} holds for all pair of vertices of any triangle of $\mathcal{T}_l$, as well as for any pair of vertices of a couple of adjacent triangles of $\mathcal{T}_l$. We now set $(\tilde{w}_\delta,\tilde{v}_\delta)(z)=(w_\delta, v_\delta)(z)$ for any vertex $z$ of a triangle in $\mathcal{T}_l$. We will claim the following:
    
    \begin{prop}\label{triangulationprop}
    		Let $\Omega^*$, $\Omega$, $l$, $(w_\delta,v_\delta)$ be as above, then there exist some universal constant $\eta>0$ and a triangulation $\tilde{\mathcal{T}}_l$ of $\Omega^*$ of minimum side $l_{\min{}}\ge \eta l$ and some $(\tilde{w}_\delta,\tilde{v}_\delta)(z_i)$ defined on each vertex $z_i$ of $\tilde{\mathcal{T}}_l$ such that, for every $z_i, z_j$ vertices of some triangle $T\in \tilde{\mathcal{T}}_l$ with $T\cap\Omega\neq \emptyset$:
    \begin{gather}
        \label{finalvonkarman}\frac{1}{2} \abs{\tilde{v}_\delta(z_i)-\tilde{v}_\delta(z_j)}^2+(\tilde{w}_\delta(z_i)-\tilde{w}_\delta(z_j))\cdot(z_i-z_j) =0\\
        \label{finalestimate}\lVert (w_\delta, v_\delta)(z_i)-(\tilde{w}_\delta, \tilde{v}_\delta)(z_i)\rVert\le C_1\delta\\
        \label{finaljump}  \abs{\tilde{v}_\delta(z_i)-\tilde{v}_\delta(z_j)}\le C_2 l\\
         \label{univbound} \Gamma(\tilde{\mathcal{T}}_l)\le C_3
    \end{gather}
     for some $C_1$ and $C_2$ depending only on the $\mathcal{C}^{1}$ norm of $(w_\delta, v_\delta)$ for $l$ small enough, and thus only on $\norm{(w,v)}_{\mathcal{C}^{1}(\Omega^*)}$, and some universal constant $C_3>0$.
    \end{prop}

 If such a construction exists (which we will proceed to show in the next section), then we may construct by affine extension on each triangle of $\tilde{\mathcal{T}}_l$ a map $(\tilde{w}_\delta,\tilde{v}_\delta)$ that is piecewise affine and such that:
    \begin{equation*}
        \frac{1}{2}\nabla \tilde{v}_\delta \otimes \nabla \tilde{v}_\delta+\sym \nabla \tilde{w}_\delta =0
    \end{equation*}
    where we have used \eqref{finalvonkarman} and Lemma \ref{interpolation} followed by the third point of Lemma \ref{preplemma2}.
    
    We can choose $l$ small enough so that from \eqref{finaljump} we get the estimates on the gradient of $\tilde{v}_\delta$ and the bound on its maximum jump; also from \eqref{univbound} we get the desired bound on the degeneracy factor. In fact, as long as the minimum side of the triangulation has length $\eta l$ for some universal $\eta $ we may see that there exists some constant $C=C_\epsilon(\norm{(w,v)}_{\mathcal{C}^1})$ such that $\abs{\nabla \tilde{v}_\delta} \le C $. We also get that :
\[
    \lambda_\varepsilon\le 2\|\nabla\tilde v_\varepsilon\|_{L^\infty}\le 2C_\varepsilon(\|(w_0,v_0)\|_{C^1}),
\]
    for the norm of the maximum jump $\lambda_\epsilon$ of $\nabla \tilde{v}_{\delta_\epsilon}$, proving the desired bound. We have also proven the gradient bound $\norm{\nabla\tilde{v}_\delta}_{L^\infty}\le C_ \epsilon(\norm{(w,v)}_{\mathcal{C}^1})$.
    
    Finally, from \eqref{finalestimate}, we would have:
    \[
        \lVert (w_\delta,v_\delta)-(\tilde{w}_\delta,\tilde{v}_\delta)\rVert_{L^\infty(\Omega ; \mathbb{R})}\le C\delta+\norm{\nabla v_\delta}l\le C_{\norm{(w,v)}_{\mathcal{C}^{1}(\Omega^*)}} \delta
    \]
  And ultimately we get:
    \begin{equation*}
    \begin{split}
        \lVert (w,v)-(\tilde{w}_\delta,\tilde{v}_\delta)\rVert_{L^\infty(\Omega ; \mathbb{R})}\le \lVert (w,v)-(w_\delta,v_\delta)\rVert_{L^\infty(\Omega ; \mathbb{R})}+\lVert (w_\delta,v_\delta)-(\tilde{w}_\delta,\tilde{v}_\delta)\rVert_{L^\infty(\Omega ; \mathbb{R})}\le \\
        \le C_{\norm{(w,v)}_{\mathcal{C}^1(\Omega^*)}}\delta \diam{\Omega}+C_{\norm{(w,v)}_{\mathcal{C}^1(\Omega^*)}}\delta\le C\delta
    \end{split}
    \end{equation*}
    where $C$ only depends on $(w,v)$ and $\Omega^*$. Hence, since $\delta$ is fully arbitrary and has been defined after the choice of $(w,v)$, we have that the thesis holds for $\delta$ small enough on $\Omega$.
\end{proof}

\subsection{Construction of the piecewise affine map}

We will now prove Proposition \ref{triangulationprop} and complete the proof of Theorem \ref{mainthm3} by showing our construction on a single fold $[abc]\cup[acd]$. 

Let us consider some $(w,v)\in \mathcal{C}^{1,\beta}(\Omega;\mathbb{R}^2\times\mathbb{R})$ satisfying:
\begin{equation}\label{deltacond}
    \frac{1}{2} \nabla v \otimes \nabla v+\sym \nabla w = -(2-\delta)\delta \id_2 
\end{equation}

As mentioned in the above section, it suffices to find a sub-triangulation of $\mathcal{T}_l$ and some map $(\tilde{w},\tilde{v})$ such that $\eqref{finalvonkarman}$ is satisfied on each triangle. We will proceed via a geometrical approach.

\begin{figure}[t]
\centering
\includegraphics[width=10em]{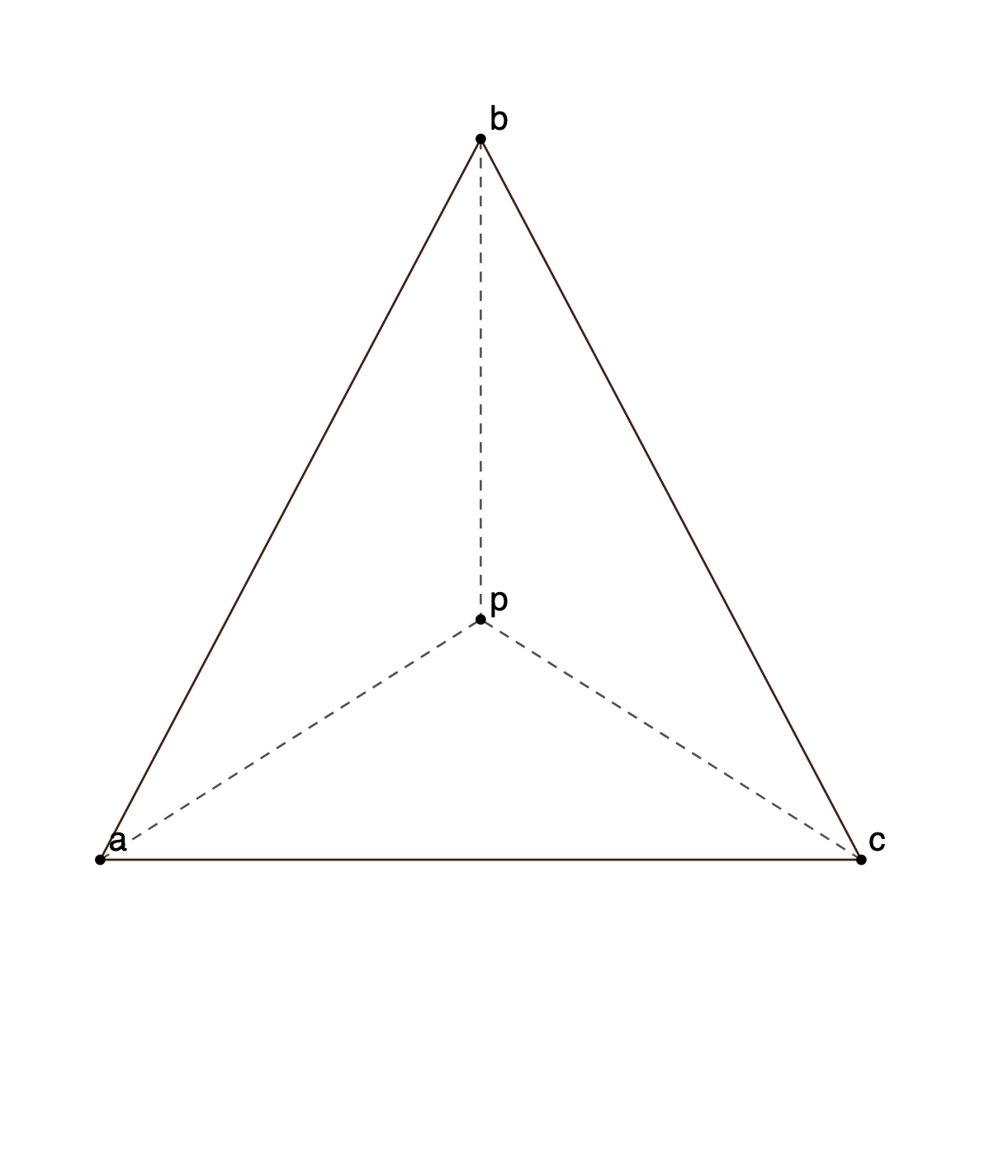}
\includegraphics[width=15em]{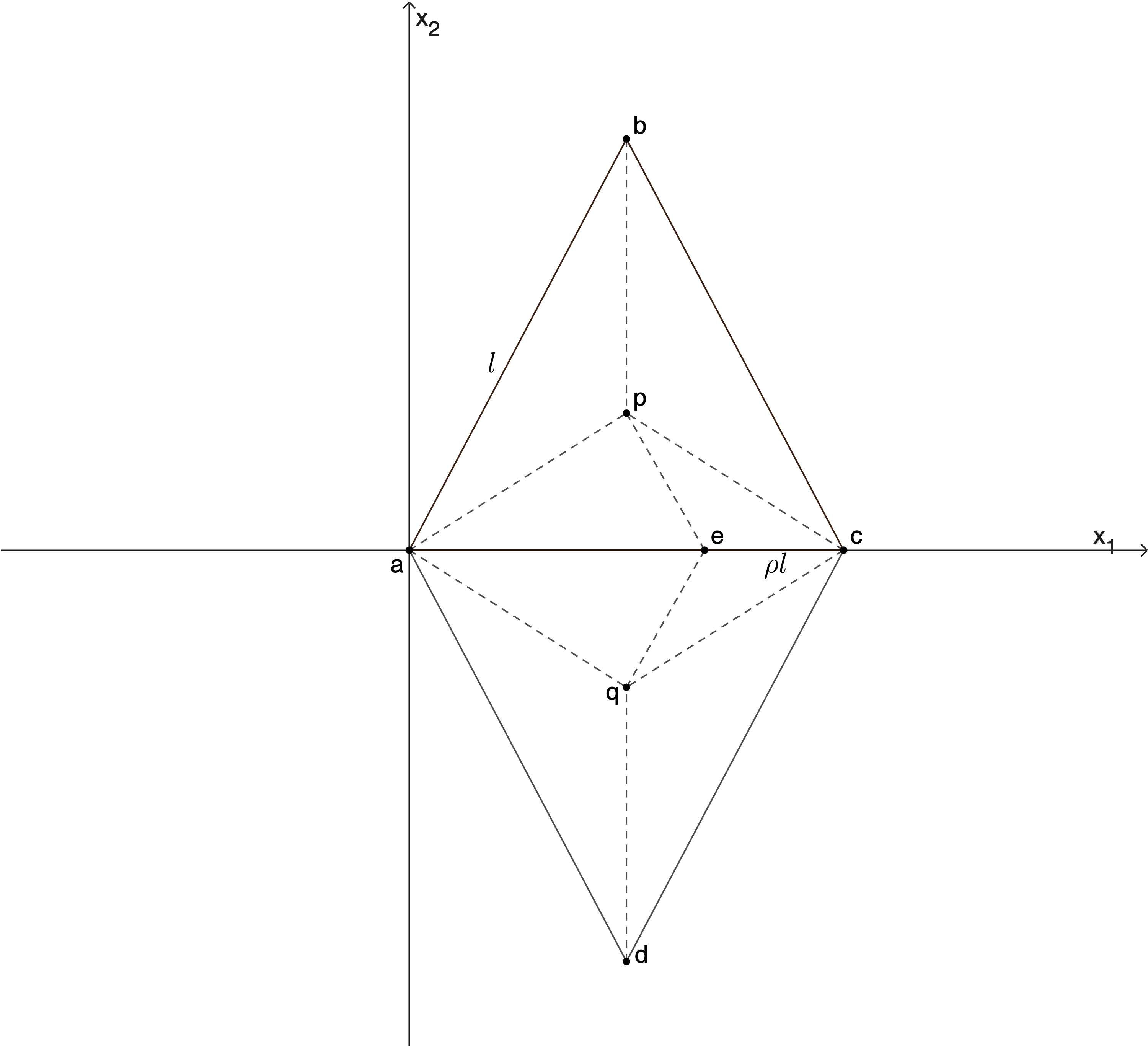}
\caption{\label{figure4}A setup of the model for the domain of the single fold considered.}
\end{figure}

Given an equilateral triangle $[abc]\in \mathcal{T}_l$ we will start by taking its barycenter $p$ and choose $(\tilde{w},\tilde{v})(p)$ in such a way that $(\tilde{w},\tilde{v})$ satisfies $\eqref{cond}$ on the segments $[ap],[bp]$ and $[cp]$. Then, we proceed to choose a point on each side of $[abc]$ and fix $(\tilde{w},\tilde{v})$ accordingly on that point. In fact, if we were to consider a pair of adjacent equilateral triangles $[abc]$ and $[acd]$ with barycenters $p$ and $q$ respectively, then, to provide a suitable subtriangulation of $[apcq]$, it suffices to choose $e\in[a,c]$ and $(\tilde{w},\tilde{v})(e)$ in such a way that \eqref{cond} is satisfied on the end points of the segments $[ae],[ce],[pe]$ and $[qe]$ and hence on the sides of the triangles $[ape],[cpe],[aqe]$ and $[cqe]$. Then, by repeating the same construction on each common side of each pair of adjacent triangles in $\mathcal{T}_l$ we get a new triangulation $\tilde{\mathcal{T}}_l$ (by taking the triangles $[ape],[cpe],[aqe]$ and $[cqe]$, see Figure \ref{figure4}), and by affine interpolation on each triangle, by Lemma \ref{interpolation} we would get $(\tilde{w},\tilde{v})$ satisfying $\eqref{cond}$ inside each triangle as well. Since we are each time applying modifications on inner points of regions with disjoint interiors, the iterations of our construction are indeed compatible with each other. Finally, for triangles at the border of $\Omega_{\mathcal{T}_l}$, we may see how those triangles have empty intersection with $\Omega$ for our assumptions on $l$. Thus if for the triangle $[abc]$ with barycenter $p$ the edge $[ab]$ is part of $\partial\Omega_{\mathcal{T}_l}$, we will just take the triangle $[abp]$ in $\tilde{\mathcal{T}}_l$, since we do not need it to satisfy any particular constraint.

So let us show that such a construction effectively exists on the triangles $[abc]$ and $[acd]$ and that it satisfies the bounds of Proposition \ref{triangulationprop}. For the first part let us take a triangle $[a_1a_2a_3]$. We start with $(\tilde{w},\tilde{v})(a_i)=(w,v)(a_i)$, for $i=1,2,3$.

By \eqref{delta3} on each pair of vertices of the triangle, we have that:
\begin{gather}\label{ab}
    -2\kappa(2-\delta)\delta l^2\le \abs{v(a_i)-v(a_j)}^2+2(w(a_i)-w(a_j))\cdot(a_i-a_j) \le -2(2-\delta)\delta l^2
\end{gather}

for $i,j\in\{1,2,3\}; i\neq j$. Let $p$ be the barycenter of the triangle, then we would like to choose $(\tilde{w},\tilde{v})(p)$ such that:
\begin{gather}\label{ap}
    \abs{\tilde{v}(p)-v(a_i)}^2+2(p-a_i)\cdot (\tilde{w}(p)-w(a_i))=0
\end{gather}
for $i=1,2,3$. By summing up \eqref{ap} for all values of $i$ and remembering that $p=\frac{a_1+a_2+a_3}{3}$ we get:
\begin{equation*}
     \sum_{i=1}^3 \abs{\tilde{v}(p)-v(a_i)}^2-2(p-a_i)\cdot w(a_i)=0
\end{equation*}
which solved for $\tilde{v}(p)$ gives us:
\begin{equation}\label{vp1}
    \tilde{v}(p)=\frac{1}{3}\left[\left(\sum_{i=1}^3 v(a_i)\right)\pm \sqrt{\sum_{\substack{i,j=1\\i\neq j}}^3 v(a_i)(v(a_j)-v(a_i))+3(p-a_i)\cdot w(a_i)}\right]
\end{equation}
where the cyclical sums are performed over $a_1, a_2$ and $a_3$. We may also notice that $\bar{v}(p)\coloneqq\frac{1}{3}\left(\sum_{i=1}^3 v(a_i)\right)$ is the affine interpolation of $v$ over the triangle $[a_1a_2a_3]$ evaluated in $p$. For reasons that will be clear later (see the inequalities leading to \eqref{pq2}), it will be necessary in our following construction to displace $\tilde{v}$ from the affine interpolated position in the same orientation for all triangles in $\mathcal{T}_l$, specifically (for simplicity) we will choose the positive sign in our computations for every barycenter. 

Finally, to prove that $\tilde{v}(p)$ exists we may sum up the right-hand side inequalities of \eqref{ab} for all distinct pair of vertices, and via some straightforward algebraic manipulations, we get:
\[
    D\coloneqq \sum_{\substack{i,j=1\\i\neq j}}^3 v(a_i)(v(a_j)-v(a_i))+3(p-a_i)\cdot w(a_i)\ge 6(2-\delta)\delta l^2> 0
\]
And we have $\tilde{v}(p)=\bar{v}(p)+ \frac{1}{3}\sqrt{D}$. We also, get, from the left-hand side inequalities in \eqref{ab} for all $i\neq j$ that:
\begin{equation}\label{bigD}
    D\le 6\kappa (2-\delta)\delta l^2
\end{equation}
Since the vectors $p-a_1$ and $p-a_2$ are linearly independent there exists exactly one solution to the system derived by \eqref{ap} for both $\tilde{w}_1(p)$ and $\tilde{w}_2(p)$, thus $\tilde{w}(p)=(\tilde{w}_1(p),\tilde{w}_2(p))$ can be then directly deduced by substituting and performing straightforward computations.

Now let us consider the adjacent triangles $[abc]$ and $[acd]$. We may define on their respective barycenters $p$ and $q$ the maps $(\tilde{w},\tilde{v})(p), (\tilde{w},\tilde{v})(q)$ as done above, and now work on the triangles $[apc],[aqc]$ to fix a point $e$ on $[ac]$ as explained before. To do so, we will work in cartesian coordinates for a more direct approach. This will require however a series of simplifications which we will attain through a series of transformations that can be seen more in detail in the Appendix. Of course, when applying the inverse transformations to go back to the original setup, we want to ensure that distances between the functions and their jumps are modified in a way that is consistent with \eqref{finalestimate} and \eqref{finaljump} and also that \eqref{deltacond} remains unchanged.

Firstly, without any loss of consistency we may apply a translation to the domain, so that $a=(0,0)$, and another translation to the image, so that $(w,v)(a)=(0,0,0)$. Now as seen in Lemma \ref{FR} we may apply any rotation $R$ to the domain via the transformation $\mathcal{D}_R$, while leaving \eqref{deltacond} unchanged (as we can see from \eqref{FR1bis}). As we go back to the original domain we will apply $\mathcal{D}_{R^{-1}}$ as seen in \eqref{FR3}, thus by \eqref{FR2bis} distances between functions will not increase and any estimate on the jumps will remain consistent as well. Thus, for simplicity, we may directly assume to have the domain rotated in such a way that $c=(l,0)$, $b=(l/2,l\sqrt{3}/2)$ and $d=(l/2,-l\sqrt{3}/2)$.

Lastly, we may apply a transformation $\mathcal{S}_\alpha$ as seen in Lemma \ref{talpha} such that $v_{(\alpha)}(c)=0$; specifically for that to happen we need $\alpha_1=-v(c)/l$ and we may fix any value for $\alpha_2$, e.g. $\alpha_2=0$. We can see that this transformation preserves \eqref{deltacond}. Also, when going back to the original system through $\mathcal{S}_{-\alpha}$ (see \eqref{talpha4}) the transformation does not modify the value of $(w,v)(a)$ we just fixed and if $\norm{(\tilde{w}_{(\alpha)},\tilde{v}_{(\alpha)})-(w_{(\alpha)},v_{(\alpha)})}_{L^\infty}<C\delta$ then $\norm{(\tilde{w},\tilde{v})-(w,v)}_{L^\infty}<C(1+\norm{\alpha})\delta$ from \eqref{talpha2} and \eqref{talpha3}. Since we have that via a Taylor expansion 
\[
    v(c)=v(a)+\nabla v(a)(c-a)+O(\abs{c-a}^2)=\nabla v(a)(c-a)+O(\abs{c-a}^{1+\beta})
\]
Then it holds:
\[
    \abs{\alpha_1}=\frac{\abs{v(c)}}{l}\le \norm{\nabla v}_{L^\infty}+O(l^\beta)
\]
So finally for $l$ small enough:
\[
    \norm{(\tilde{w},\tilde{v})-(w,v)}_{L^\infty}<C(2+\norm{\nabla v}_{L^\infty})\delta
\]
which satisfies the requirements we wanted for \eqref{finalestimate}. In the same way from \eqref{talphagrad} and \eqref{talphagrad2} we may see that if $\abs{(\tilde{w},\tilde{v})(x)-(\tilde{w},\tilde{v})(y)}\le C l$ for any two adjacent vertices $x,y$, then:
\[
    \abs{(\tilde{w}_{(\alpha)},\tilde{v}_{(\alpha)})(x)-(\tilde{w}_{(\alpha)},\tilde{v}_{(\alpha)})(y)}<C(2+\norm{\nabla v}_{L^\infty}) l
\]
 which gives us consistency for \eqref{finaljump} as well. Hence we may assume directly that $v(c)=0$.


With these new data we get, by substituting in \eqref{vp1} and \eqref{ap} for $a_1, a_2, a_3 = a, b, c$, with $a+b+c=3p$, and performing computations: 
\begin{gather*}
    \tilde{v}(p)=\frac{v(b)+\sqrt{D_p}}{3}, \quad \quad \quad \quad \quad
    \tilde{w}_1(p)=\frac{3w_1(c)-\sqrt{3}w_2(c)}{6},\\
    \tilde{w}_2(p)=\frac{2\sqrt{3}v(b)^2-3\sqrt{3}lw_1(c)+12l w_2(b)+3lw_2(c)-4\sqrt{3}v(b)\sqrt{D_p}}{18l},
\end{gather*}
with:
\begin{equation}\label{dp1}
	D_p=-2v(b)^2-2(b-c)\cdot (w(b)-w(c))-2b\cdot w(b)-2c\cdot w(c)\ge 4(2-\delta)\delta l^2-2c\cdot w(c)
\end{equation}

Also we may see that, by using again the Taylor expansion for $(w,v)$ with respect to the point $a=0$, for $l$ small enough:
\begin{equation}\label{estimate1}
\begin{split}
    \abs{(\tilde{w},\tilde{v})(p)-(w,v)(p)} &\le \abs{(\tilde{w},\tilde{v})(p)-(w,v)(a)}+\abs{(w,v)(a)-(w,v)(p)}\\
    &\le C_{\norm{(\nabla w, \nabla v)}_\infty} l+\norm{( w, v)}_{\mathcal{C}^{1,\beta}}O(l^{1+\beta})\\&\le C_{\norm{( w, v)}_{\mathcal{C}^{1}}} \delta
\end{split}
\end{equation}

On the other hand on the triangle $[acd]$ we get:
\[
	\tilde{v}(q)=\frac{v(d)+\sqrt{D_q}}{3}
\]
with:
\begin{equation}\label{dq1}
	D_q=-2v(d)^2-2(d-c)\cdot (w(d)-w(c))-2b\cdot w(d)-2c\cdot w(c)\ge 4(2-\delta)\delta l^2-2c\cdot w(c)
\end{equation}
And a similar estimate as \eqref{estimate1} holds for $q$ as well. So, from \eqref{dp1} and \eqref{dq1}:
\begin{equation}\label{vpvq1}
	\sqrt{D_p}+\sqrt{D_q}\ge 2\sqrt{4(2-\delta)\delta l^2-2c\cdot w(c)}
\end{equation}


We also have by using Lemma \ref{preplemma2} (i) on $[ac]$ and from our construction of $p, q$ on the sides $[ap],[pc],[aq]$ and $[qc]$ respectively:
\begin{gather}\label{ac2}
     2w(c)\cdot c\le -2(2-\delta)\delta l^2\\
    \label{ap2} \tilde{v}(p)^2+2\tilde{w}(p)\cdot  p=0\\
    \label{pc2} \tilde{v}(p)^2+2(\tilde{w}(p)-w(c))\cdot(p-c)=0\\
    \label{aq2} \tilde{v}(q)^2+2\tilde{w}(q)\cdot  q=0\\
    \label{qc2} \tilde{v}(q)^2+2(\tilde{w}(q)-w(c))\cdot(q-c)=0
\end{gather}
We may see right away that subtracting \eqref{ap2} from \eqref{pc2} and \eqref{aq2} from \eqref{qc2} and summing up the two resulting equations we get:
\begin{equation}\label{pqc}
	(\tilde{w}(p)+\tilde{w}(q))\cdot c = w(c)\cdot c
\end{equation}
where we have used $p+q=c$ (from the symmetry of our construction). It follows that:
\begin{equation}\label{trick}
    \tilde{w}_1(p)+\tilde{w}_1(q)=w_1(c)
\end{equation}

Lastly, we notice that by applying Lemma \ref{preplemma2} (ii) on $[ac]$ we get:
\begin{equation}\label{ac3}
	2w(c)\cdot c\ge -2\kappa (2-\delta)\delta l^2
\end{equation}
and by \eqref{bigD}:
\begin{equation}\label{doubleD}
	\sqrt{D_p}+\sqrt{D_q}\le 2 l\sqrt{6\kappa(2-\delta)\delta }
\end{equation}
also, by symmetry of our construction:
\begin{equation}\label{smalls}
	\abs{\frac{v(b)+v(d)}{3}}=\abs{\frac{v(b)-v(a)+v(d)-v(c)}{3}}\le \norm{v}_{\mathcal{C}^{1,\beta}} l^{1+\beta}
\end{equation}

Before proceeding with our edge construction we will perform some further algebraic manipulation to obtain some bound on:
\[
	\Delta(p,q)=(\tilde{v}(p)-\tilde{v}(q))^2+2(\tilde{w}(p)-\tilde{w}(q))\cdot(p-q)
\]
First of all we see, by applying \eqref{ap2} and \eqref{aq2} followed by \eqref{pqc}:
\begin{equation*}
	\begin{split}
		\Delta(p,q)&=-(\tilde{v}(p)+\tilde{v}(q))^2+2\tilde{v}(p)^2+2\tilde{v}(q)^2+2\tilde{w}(p)\cdot p +2\tilde{w}(q)\cdot q -2\tilde{w}(p)\cdot q-2\tilde{w}(q)\cdot p\\
		&=-(\tilde{v}(p)+\tilde{v}(q))^2-2\tilde{w}(p)\cdot p -2\tilde{w}(q)\cdot q -2\tilde{w}(p)\cdot q-2\tilde{w}(q)\cdot p\\
		&=-(\tilde{v}(p)+\tilde{v}(q))^2-2(\tilde{w}(p)+\tilde{w}(q))\cdot c\\
		&=-(\tilde{v}(p)+\tilde{v}(q))^2-2w(c)\cdot c
	\end{split}
\end{equation*}
Hence by applying \eqref{vpvq1} and \eqref{smalls}, followed by \eqref{ac3} and \eqref{doubleD}, by remembering that we chose $\kappa < \frac{3}{2}$:
\[
	\begin{split}
		-\Delta(p,q)&=\left(\frac{v(b)+v(d)+\sqrt{D_p}+\sqrt{D_q}}{3}\right)^2+2w(c)\cdot c\\
		&\ge -\frac{2}{9}\abs{v(b)+v(d)}(\sqrt{D_p}+\sqrt{D_q})+\frac{4}{9}\left(\sqrt{4(2-\delta)\delta l^2-2c\cdot w(c)}\right)^2+2w(c)\cdot c\\
		&\ge -\frac{2}{3}(\sqrt{D_p}+\sqrt{D_q})\norm{v}_{\mathcal{C}^{1,\beta}} l^{1+\beta}+\frac{16}{9}(2-\delta)\delta l^2+\frac{10}{9}c\cdot w(c)\\
		&\ge -\frac{4}{3}\sqrt{6\kappa(2-\delta)\delta }\norm{v}_{\mathcal{C}^{1,\beta}} l^{2+\beta}+\left(\frac{16}{9}-\frac{10}{9} \kappa\right)(2-\delta)\delta l^2\\
		&\ge -4\sqrt{(2-\delta)\delta }\norm{v}_{\mathcal{C}^{1,\beta}} l^{2+\beta}+\frac{1}{9}(2-\delta)\delta l^2\ge \frac{1}{18}(2-\delta)\delta l^2>0
	\end{split}
\]
as long as we choose $l$ small enough, specifically such that $\norm{v}_{\mathcal{C}^{1,\beta}} l^{\beta}\le \frac{1}{72}(2-\delta)\delta$. Hence we ultimately get:
\begin{equation}\label{pq2}
    \Delta(p,q)=(\tilde{v}(p)-\tilde{v}(q))^2+2(\tilde{w}(p)-\tilde{w}(q))\cdot(p-q)< 0
\end{equation}

Let us now proceed with our construction. We would like to find some point $e\in [ac]$ and some values $(\tilde{w},\tilde{v})(e)$ such that:
\begin{gather}
    \label{ae} \tilde{v}(e)^2+2\tilde{w}(e)\cdot  e=0\\
    \label{ec} \tilde{v}(e)^2+2(\tilde{w}(e)-w(c))\cdot(e-c)=0\\
    \label{pe} (\tilde{v}(p)-\tilde{v}(e))^2+2(\tilde{w}(p)-\tilde{w}(e))\cdot  (p-e)=0\\
    \label{qe} (\tilde{v}(q)-\tilde{v}(e))^2+2(\tilde{w}(q)-\tilde{w}(e))\cdot(q-e)=0
\end{gather}
We may write $e=(1-\rho)a+\rho c=\rho c$ and to prove its existence it suffices to have $\rho\in (0,1)$. So we get from \eqref{ae} and \eqref{ec}:
\begin{gather*}
    \tilde{v}(e)^2+2\rho\tilde{w}(e)\cdot  c=0\\
    \tilde{v}(e)^2+2(\rho-1)(\tilde{w}(e)-w(c))\cdot c=0
\end{gather*}
By multiplying the former by $(1-\rho)$ and the latter by $\rho$ and summing them up we have:
\begin{equation}\label{vrho}
    \tilde{v}(e)^2+2\rho(1-\rho)w(c)\cdot c=0
\end{equation}
Hence:
\[
    \tilde{v}(e)=\pm\sqrt{-2\rho(1-\rho)w(c)\cdot c}
\]
which exists as long as the radicand is positive, hence by \eqref{ac2} for its existence it suffices to prove that $\rho\in(0,1)$. We may also estimate:
\begin{equation}\label{estimate2}
    \abs{\tilde{v}(e)-v(e)}\le \sqrt{l\abs{w(c)}} + \abs{v(e)} \le C_{\norm{ (w,v)}_{\mathcal{C}^{1}}} l
\end{equation}
On the other hand, by subtracting \eqref{ae} from \eqref{ec}, \eqref{pe} and \eqref{qe} we get:
\begin{gather}\label{wrho}
    \tilde{w}(e)\cdot c=(1-\rho)w(c)\cdot c\\
    \label{pe2} \tilde{v}(p)^2-2\tilde{v}(p)\tilde{v}(e)+2\tilde{w}(p)\cdot  (p-\rho c) -2\tilde{w}(e)\cdot p=0\\
    \label{qe2} \tilde{v}(q)^2-2\tilde{v}(q)\tilde{v}(e)+2\tilde{w}(q)\cdot  (q-\rho c) -2\tilde{w}(e)\cdot q=0
\end{gather}
We may simplify \eqref{pe2} and \eqref{qe2} by using \eqref{ap2} and \eqref{aq2} respectively and get:
\begin{gather}
    \label{pe3} \tilde{v}(p)\tilde{v}(e)+\rho\tilde{w}(p)\cdot  c +\tilde{w}(e)\cdot p=0\\
    \label{qe3} \tilde{v}(q)\tilde{v}(e)+\rho\tilde{w}(q)\cdot  c +\tilde{w}(e)\cdot q=0
\end{gather}
We have that $p+q=c$, also by \eqref{trick} we get $(\tilde{w}(p)+\tilde{w}(q))\cdot  c=w(c)\cdot c$ . Consequently, by summing up \eqref{pe3} and \eqref{qe3}:
\[
    (\tilde{v}(p)+\tilde{v}(q))\tilde{v}(e)+\rho w(c)\cdot  c +\tilde{w}(e)\cdot c=0
\]
By substituting the expression we found for $\tilde{w}(e)\cdot c$ in \eqref{wrho} and isolating $\tilde{v}(e)$ we finally have:
\begin{equation}\label{vefinal}
    \tilde{v}(e)= -\frac{w(c)\cdot c}{\tilde{v}(p)+\tilde{v}(q)}
\end{equation}
as long as $\tilde{v}(p)+\tilde{v}(q)\neq0$. We will now substitute this newfound value for $\tilde{v}(e)$ in \eqref{vrho} to finally compute $\rho$. We ultimately get:
\[
    \rho= \frac{1}{2}\pm \frac{\sqrt{(\tilde{v}(p)+\tilde{v}(q))^2+2w(c)\cdot c}}{2\abs{\tilde{v}(p)+\tilde{v}(q)}}
\]
We may see that it is symmetrical around the ideal value $\rho=\frac{1}{2}$. Furthermore, we may see that, by adding \eqref{pq2} to the radicand and using \eqref{ap2} and \eqref{aq2}:
\[
    \begin{split}
        (\tilde{v}(p)+\tilde{v}(q))^2+2w(c)\cdot c &=(\tilde{v}(p)+\tilde{v}(q))^2+2(\tilde{w}(p)+\tilde{w}(q))\cdot (p+q)\\
        &\ge 2\left[\tilde{v}(p)^2+\tilde{v}(q)^2+(\tilde{w}(p)+\tilde{w}(q))\cdot (p+q)+(\tilde{w}(p)-\tilde{w}(q))\cdot(p-q)\right]\\
        &=2\tilde{v}(p)^2+4\tilde{w}(p)\cdot p+2\tilde{v}(q)^2+4\tilde{w}(q)\cdot q=0
    \end{split}
\]
so the radicand is positive. Furthermore this proves that, by using \eqref{ac2}:
\[
    (\tilde{v}(p)+\tilde{v}(q))^2\ge -2w(c)\cdot c \ge 2(2-\delta)\delta l^2
\]
Hence, $\abs{\tilde{v}(p)+\tilde{v}(q)}\ge C_1 l\sqrt{\delta}>0$ for some universal constant $C_1$. This proves that we already had a priori $\tilde{v}(p)+\tilde{v}(q)\neq 0$ and hence $\rho$ exists. More precisely, we may see that, since $2w(c)\cdot c<-C_2\delta l^2$ for some universal constant $C_2>0$:
\[
    \frac{\sqrt{(\tilde{v}(p)+\tilde{v}(q))^2+2w(c)\cdot c}}{\abs{\tilde{v}(p)+\tilde{v}(q)}}\le \frac{\sqrt{(\tilde{v}(p)+\tilde{v}(q))^2-C_2\delta l^2}}{\abs{\tilde{v}(p)+\tilde{v}(q)}}\le 1
\]
hence $\rho\in (0,1)$ and $\tilde{v}(e)$ exists. 


We would now like to put some bounds on the minimum side of the triangles constructed and prove that \eqref{univbound} holds. Since the maximum sides of the triangles defined and $\nu=\nu(\tilde{\mathcal{T}}_l)$ are bounded for our construction, for the degeneracy factor it suffices to show that $\frac{1}{\tau_{\min{}}}<C$, which is true as long as there exists some universal $\eta>0$ such that $\eta<\rho<1-\eta$: in fact that would prove that $\abs{[ae]},\abs{[ce]}>\eta l$ and bound the minimum side of the constructed triangles, and thus, since one of their angles is fixed, their minimum slope as well. To see that, we will look at the expression of $\rho$. We have:
\[
    \rho=\frac{1}{2}\pm \frac{1}{2} \sqrt{1+\frac{2w(c)\cdot c}{(\tilde{v}(p)+\tilde{v}(q))^2}}
\]
The problem arises if the radicand is arbitrarily close to 1. Hence it suffices to prove that there exists some universal $k>0$ such that:
\[
    -\frac{w(c)\cdot c}{(\tilde{v}(p)+\tilde{v}(q))^2}\ge k
\]
We have that, by \eqref{ac2}, $-w(c)\cdot c\ge (2-\delta)\delta l^2$. On the other hand we can see that, by \eqref{vp1} and \eqref{bigD}:
\[
    (\tilde{v}(p)+\tilde{v}(q))^2\le 2\left[(v(b)+v(d))^2+12\kappa(2-\delta)\delta l^2\right]
\]
with $\kappa\in(1,3/2)$. Now, since $v(a)=v(c)=0$, there exists some point $e^*\in [ac]$ such that $\nabla v(e^*)\cdot c =0$. Let us look at $v(b)$ and $v(d)$ via a Taylor expansion with respect to $a$ and get, since $v\in \mathcal{C}^{1,\beta}$:
\[
    (v(b)+v(d))^2=\left(\nabla v(a) \cdot (b+d)+O\left(l^{1+\beta}\right)\right)^2\le C_1\left(\abs{e^*}^\beta 
l+O\left(l^{1+\beta}\right)\right)^2\le C_1l^{2+2\beta}
\]
For some constant $C_1$ depending on the $\mathcal{C}^{1,\beta}$ norm of $v$. By choosing $l$ small enough we have $C_1 l^{2\beta}<\delta$ hence for at most some some universal constant $C$:
\[
    (\tilde{v}(p)+\tilde{v}(q))^2\le C\delta l^2
\]
and we get in the end:
\begin{equation}\label{finalunivestimate}
    -\frac{w(c)\cdot c}{(\tilde{v}(p)+\tilde{v}(q))^2}\ge \frac{(2-\delta)\delta l^2}{C\delta l^2}\ge k>0
\end{equation}
for some universal $k$, which in turn proves \eqref{univbound}.


Finally we may now look at the simplified systems of equations equivalent to \eqref{ae}-\eqref{qe}:
\[
\begin{cases}
\tilde{v}(e)^2+2\tilde{w}(e)\cdot  e=0\\
\tilde{w}(e)\cdot c=(1-\rho)w(c)\cdot c\\
\tilde{v}(p)\tilde{v}(e)+\rho\tilde{w}(p)\cdot  c +\tilde{w}(e)\cdot p=0\\
\tilde{v}(q)\tilde{v}(e)+\rho\tilde{w}(q)\cdot  c +\tilde{w}(e)\cdot q=0
\end{cases}
\]
By substitution in the equations above we find the values for $\tilde{w}_1(e)$ and $\tilde{w}_2(e)$ which at this point fully show the existence of our construction. Similar bounds as \eqref{estimate2} can be found by homogeneity of the equations taken into consideration, and ultimately we get:
\[
    \abs{(\tilde{w},\tilde{v})(e)-(w,v)(e)}\le C_{\norm{ (w,v)}_{\mathcal{C}^{1}}}l
\]
which for $l$ small enough gives us \eqref{finalestimate} as we wanted.

For the maximum jump of the constructed $\tilde{v}$ over the new triangulation $\tilde{\mathcal{T}}_l$, from \eqref{vp1}, by using \eqref{bigD} and since $v(b)\sim Cl$ for some $C$ depending on $\norm{(w,v)}_{\mathcal{C}^{1}}$, then we have $\abs{\tilde{v}(p)},\abs{\tilde{v}(q)}\sim C l$. Meanwhile from \eqref{estimate2} we get that $\abs{\tilde{v}(e)}\le  C l $, so, by using \eqref{finalunivestimate} to get a minimum bound of $\eta l$ for the side length, we have that \eqref{finaljump} holds.

%
%
%
%
%
%
%

%
%
%
%
%
%
%
%
%
%
%
%
%
%

\section{Approximation of a piecewise affine solution to the von Kármán constraint}

\begin{figure}[b]
\centering
\includegraphics[width=15em]{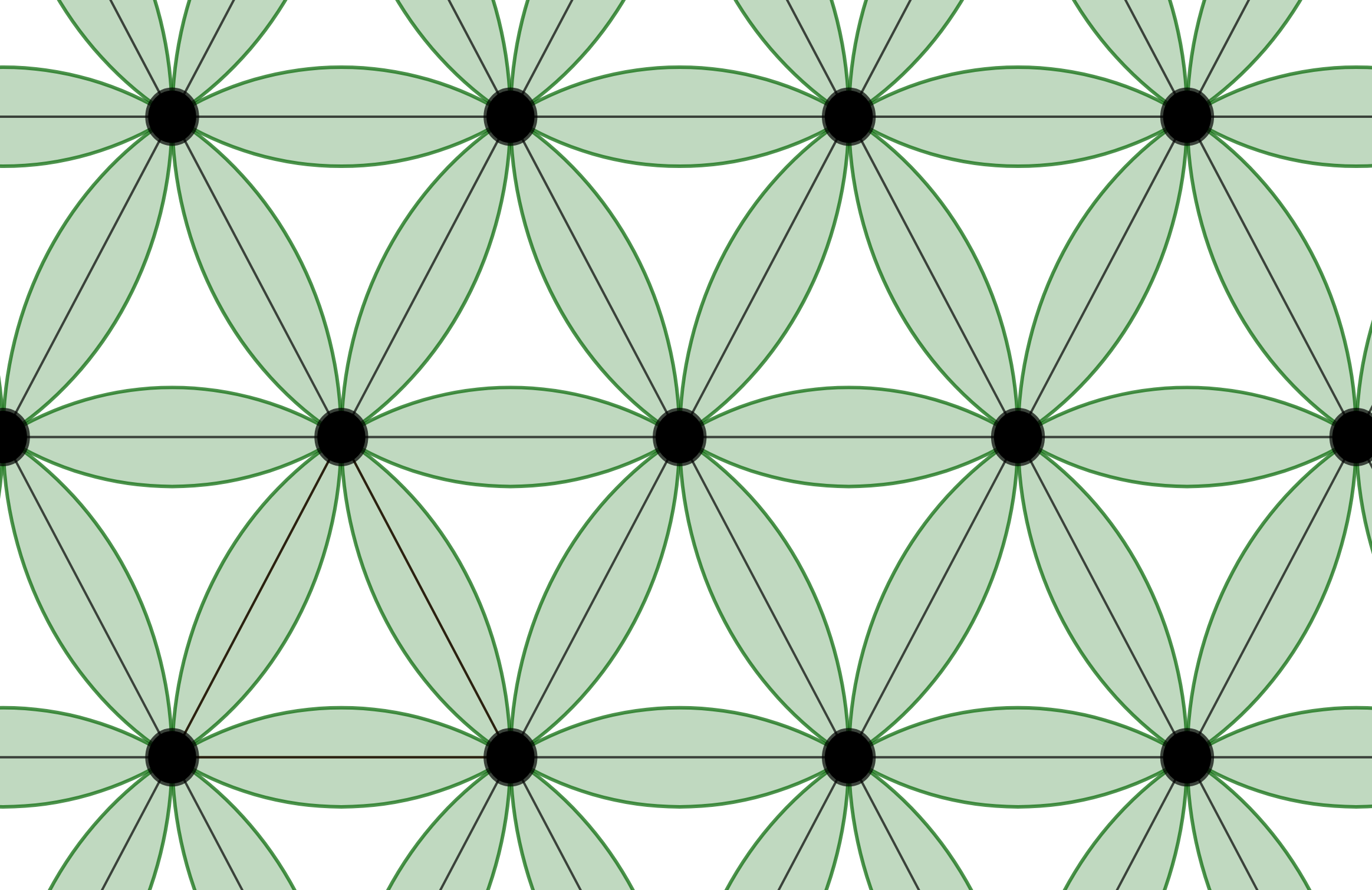}
\caption{\label{construction} This picture gives an idea of the construction: on this equilateral triangulation, the green regions represent the ones on which we would like to apply our modifications. The black circles around the vertices are the regions we will treat separately to avoid singularities.}
\end{figure}

We now focus on modifying the piecewise affine map we have constructed and provide a definitive proof to Theorem \ref{mainthm}, adapting the ridge construction used for example in \cite{contimaggi} and \cite{venkataramani} to the nonlinear vK setting. Let $\Omega\subset\mathbb{R}^2$ be an open bounded Lipschitz domain and $\mathcal{T}$ be a triangulation over said domain. We will assume that we have some piecewise affine map $(w,v) \in W^{1,\infty}(\Omega, \mathbb{R}^2\times \mathbb{R})$ such that it satisfies the vK constraint

\begin{equation}\label{vonkarman}
    \frac{1}{2}\nabla v\otimes \nabla v +\sym \nabla w = 0
\end{equation}

and that $\nabla w$ and $\nabla v$ are constant on $T_i$ for each $T_i \in \mathcal T$. We now look for some modification $w_\sigma, v_\sigma$ that can closely enough approximate our initial data. We would like to "smoothen" our deformation around each edge of $\mathcal{T}$ by taking some connecting smooth functions between the two adjacent triangles and rescaling this function along this boundary region in such a way that the transitioning region shrinks around the vertices (which decreases the surface of our modification reducing the total energy and at the same time does not let our modifications interact with each other). The resulting mapping will not be smooth around the vertices, which we will in fact isolate and work on separately. The final result we would like to obtain in this Section is the following:
\begin{thm}\label{mainthm2}
Let $(w,v)\in W^{1,\infty}(\Omega; \mathbb{R}^2\times \mathbb{R})$ be a piecewise affine solution to \ref{vonkarman} with respect to some triangulation $\mathcal{T}$ of the bounded domain $\Omega$. Let $\lambda_*\ge 0$ be the norm of the maximum jump of $\nabla v$. Then for every $h$ small enough, depending only on $\Omega$ and $\mathcal{T}$ (see \eqref{hcond}), there is a displacement map $(\bar w_h,\bar v_h)\in \mathcal{C}^2(\Omega;\mathbb{R}^2\times \mathbb{R})$ such that:
\begin{gather}
    \label{estimate3.2}\norm{\bar v_h-v}_{L^\infty(\Omega)}\le  C_{\Gamma(\mathcal{T})}\lambda_*^{2/3} h^{1/3}\\
    \label{estimate3.3}\lVert \bar w_h - w\rVert_{L^\infty(\Omega)}\le C_{\Gamma(\mathcal{T})}\lambda_*^{2/3}\norm{\nabla v}_{L^\infty} h^{1/3}\\
    \label{finalenergythm} I^h[(\bar w_h,\bar v_h),\Omega]\le C_{\Gamma(\mathcal{T})} \abs{\mathcal{T}} \lambda_*^{7/3}h^{5/3} 
\end{gather}

\end{thm}
We remark that the proof of Theorem \ref{mainthm} derives directly from combining Theorems \ref{mainthm3} and \ref{mainthm2}. Specifically we apply \ref{mainthm3} with $\epsilon = r/2$ and then we have by \eqref{estimate3.2} and \eqref{estimate3.3} that $\|(\bar w_h,\bar v_h)-(\tilde w,\tilde v)\|_\infty\le C_\Gamma\lambda_*^{2/3}(1+C_\varepsilon)h^{1/3}<r/2$ for $h^*$ small enough and finally setting $C^*=C_\Gamma|\mathcal T_\varepsilon|\lambda_*^{7/3}$. It is interesting to see that a scaling of order $5/3$ necessitates this kind of specific construction where the length of the edge deformation decreases as we move towards the vertices. We can look for example at Figure \ref{construction}. A more naive approach where we cut away "rectangles" around each side would not yield these same results providing regions with an area scaling only linearly with $h$.

To prove Theorem \ref{mainthm2}, it will be necessary to first study the nature of our construction on a \textit{single fold}, or in other words a pair of adjacent triangles $T_1$, $T_2$ in our triangulation $\mathcal{T}$, and some associated piecewise affine deformation, that can be fully defined by a pair $(w,v)$ of piecewise affine displacement functions. We will define its \emph{length} $l$ as the length of the common side of $T_1$ and $T_2$, its \emph{jump} as $\lambda=\abs{\nabla v_{|T_1}-\nabla v_{|T_2}}$ and its \emph{slope} $\tau$ as the minimum slope among the two triangles with respect to their common side. We would like to work on the common edge of those triangles and "smoothen" the fold as discussed before obtaining the following result. 
\begin{lemma}\label{singlefoldthm}
    Let $(T_1\cup T_2, (w,v))$ be a single fold, with length $l$, jump $\lambda$ and slope $\tau$; also let $x,y\in \mathbb{R}^2$ be the ends of the segment $T_1\cap T_2$. Then for every $\sigma < l/8$ there exists a mapping $(w_\sigma, v_\sigma)$ such that $(w_\sigma, v_\sigma)\in W^{1,\infty}(T_1 \cup T_2;\mathbb{R}^3)\cap \mathcal{C}^2(T_1 \cup T_2\setminus \{x,y\};\mathbb{R}^3)$, where $(w_\sigma, v_\sigma)=(w, v)$ and $(\nabla w_\sigma, \nabla v_\sigma)=(\nabla w, \nabla v)$ on $\partial (T_1\cup T_2)$, and which obeys:
    \begin{gather}
        \label{2.5}\lVert v_\sigma - v \rVert_{L^\infty(T_1 \cup T_2;\mathbb{R})} \le C_\tau \abs{\lambda} l^{2/3} \sigma^{1/3}\\
        \label{2.6}\lVert w_\sigma - w \rVert_{L^\infty(T_1 \cup T_2;\mathbb{R}^2)} \le C_\tau \abs{\lambda} \norm{\nabla v}_{L^\infty} l^{2/3} \sigma^{1/3}\\
        \label{2.7}\lVert \nabla v_\sigma - \nabla v \rVert_{L^\infty(T_1 \cup T_2;\mathbb{R}^{2})} \le C_\tau \abs{\lambda},\quad   \lVert \nabla^2 v_\sigma \rVert_{L^\infty(T_1 \cup T_2;\mathbb{R}^{2\times 2} )}\le C_\tau \frac{\abs{\lambda}}{\sigma} \\
        \label{2.8}\lVert \nabla w_\sigma - \nabla w \rVert_{L^\infty(T_1 \cup T_2;\mathbb{R}^{2\times 2})} \le C_\tau \abs{\lambda}\norm{\nabla v}_{L^\infty}
    \end{gather}
    and such that:
    \begin{equation*}
    I^h[(w_\sigma, v_\sigma), (T_1 \cup T_2)_\sigma]\le C\left[\abs{\lambda}^4\tau^5l^{1/3} \sigma^{5/3}+\abs{\lambda}^2h^2\left(\tau^3+\frac{l^{1/3}}{\tau\sigma^{1/3}}\right)\right]
    \end{equation*}
    where $(T_1 \cup T_2)_\sigma$ is $(T_1 \cup T_2) \setminus (B(x,\sigma)\cup B(y,\sigma))$. Specifically, if $\sigma=h/\abs{\lambda}$:
    \begin{equation}\label{sftresult2}
        I^h[(w_h, v_h), (T_1 \cup T_2)_h]\le C_\tau \abs{\lambda}^{7/3} l^{1/3}h^{5/3}
    \end{equation}
\end{lemma}

In sections 3.1 through 3.4 we will prove Lemma \ref{singlefoldthm}, while in the following section 3.5 we will show how that construction can be extended to the full domain and prove Theorem \ref{mainthm2}.

\subsection{Setup for our base model on a single fold}

First of all let us consider a single fold $(T_1,T_2)$ in $\mathcal{T}$. Since $(w,v)$ is affine over each of those triangles we may write:
\[
    \nabla (w,v)\vert_{T_i} = (B_i, a_i)
\]
for $i=1,2$, where $a_1,a_2\in \mathbb{R}^2$ and $B_1,B_2\in\mathbb{R}^{2\times 2}$. Hence, from \eqref{vonkarman} we get:

\begin{equation}\label{vonkarmanAB}
    \frac{1}{2}( a_1\otimes a_1 - a_2\otimes  a_2) = \sym  (B_2-B_1)
\end{equation}
We may now choose $\eta$ as a unitary normal vector to $\partial T_1 \cap \partial T_2$ in $\mathbb{R}^2$. Since $(w,v)$ is of class $W^{1,2}$ over the bounded domain $\Omega$, we have some compatibility conditions over $a_i$ and $B_i$. Specifically, as seen for example in \cite{LAB}, we have that $\rk (B_2-B_1) \le 1$ and $(B_2-B_1)\eta^\perp =0$ and $(a_2-a_1)\cdot \eta^\perp=0$. As such, we get:
\begin{equation*}
    \begin{cases}
        a_2= a_1+\lambda\eta\\
        B_2=B_1+b \otimes \eta
    \end{cases}
\end{equation*}
for some $b\in\mathbb{R}^2, \lambda\in\mathbb{R}$. Hence via some straightforward computation:
\begin{equation}\label{B}
    \sym(B_2 - B_1) = \frac{1}{2}(b\otimes \eta + \eta\otimes b )
\end{equation}
\begin{equation}\label{A}
    \begin{split}
        a_1\otimes a_1 - a_2\otimes  a_2 &= a_1\otimes a_1 - (a_1+\lambda\eta)\otimes (a_1+\lambda\eta)=\\
        &=-(a_1\otimes \lambda\eta)-(\lambda\eta\otimes a_1) - \lambda^2 \eta\otimes \eta =\\
        &=-\lambda\left(a_1+\frac{\lambda\eta}{2}\right)\otimes\eta-\eta\otimes\lambda\left(a_1+\frac{\lambda\eta}{2}\right)
    \end{split}
\end{equation}

By combining \eqref{vonkarmanAB}, \eqref{B} and \eqref{A} we get that 
\begin{equation*}
    \sym(b\otimes \eta) =-\lambda\sym \left[\left(a_1+\frac{\lambda\eta}{2}\right)\otimes\eta\right]
\end{equation*}
and since $\eta$ is nonzero by definition we necessarily have that:
\begin{equation}\label{finalB}
    b= -\lambda\left(a_1 + \frac{\lambda\eta}{2}\right)
\end{equation}

Since translations do not significantly impact the conditions for our problem, and the same holds for rotations (see Lemma \ref{FR}, \eqref{FR1bis} and \eqref{FR2bis}), we may apply a simplification analogous to the one used for the geometrical construction of the previous chapter. By at most applying an isometry, in fact, we may assume $T_1\cap T_2$ to have $(0,0)\in\mathbb{R}^2$ as one of its endpoints, $\eta= e_2$, with $e_1,  e_2$ the canonical base for $\mathbb{R}^2\supset \Omega$, or in other words so that $T_1\cap T_2 = [0, le_1]$, where $l=\abs{T_1\cap T_2}$ is the length of the fold. Also, without loss of generality, we will assume $T_2$ to be contained in the half-plane $\{x\in \mathbb{R}^2\ \colon\ x_2=x\cdot e_2 < 0\}$. 

After these transformations we have, by using \eqref{finalB} (and leaving for simplicity the names of $v$ and $w$ unchanged): 
\[
    v(x)=\begin{cases}
        a\cdot x  &\textrm{for }x\in T_1\\
        a\cdot x+\lambda x_2 \quad &\textrm{for }x\in T_2
    \end{cases}
\]
\[
    w(x)=\begin{cases}
        Bx  &\textrm{for }x\in T_1\\
        Bx-\lambda x_2a-\frac{1}{2} \lambda^2 x_2 e_2\quad &\textrm{for }x\in T_2
    \end{cases}
\]
where we have for brevity renamed $a=a_1$, $B=B_1$ and we have:
\begin{equation}\label{vonkarman2}
    \frac{1}{2}a\otimes a+\sym B = 0
\end{equation}

To conclude our setup, we may see that given a construction $( w_\sigma',v_\sigma')$ on $T_1'\subset T_1$ and $T_2'\subset T_2$, with $T_1'\cap T_2' = T_1\cap T_2$ and $w_\sigma'=w, v_\sigma'=v$ on $\partial(T_1'\cup T_2')$, then a construction $( w_\sigma,v_\sigma)$ for Lemma \ref{singlefoldthm} can be directly derived by extension with $w_\sigma=w$ and $v_\sigma = v$ outside of $T_1'\cup T_2'$. As such, increasing the slope of each side of the fold only improves the estimate, hence it suffices to prove our thesis for a fold with minimum slope $\tau=\tau_{\min}<1$. We can thus work on the restriction of our domain:
\[
    Q\coloneqq\{x\in\mathbb{R}^2 \ \colon\ 0\le x_1\le l, \abs{x_2}\le \tau\min\{x_1,l-x_1\}\}
\]
We will first construct a map $(w_\sigma,v_\sigma) \in W^{1,\infty}(T_1\cup T_2; \mathbb{R}^2\times\mathbb{R})\cap \mathcal{C}^2((0,l)\times\mathbb{R};\mathbb{R}^2\times\mathbb{R})$ such that
\begin{equation}
    (w_\sigma,v_\sigma) = (w,v) \quad \textrm{on }\ \mathbb{R}^2\setminus Q
\end{equation}

and that satisfies the estimates from Lemma \ref{singlefoldthm} on
\[
    Q_\sigma \coloneqq \{(x_1,x_2)\in Q \ :\ \min\{x_1,l-x_1\}>\sigma'\},
\]
where, since $\tau < 1$, we may assume $\sigma'=\sigma/\sqrt{2}$ so that $ (T_1 \cup T_2)_\sigma \subset ((T_1 \cup T_2)\setminus Q)\cup Q_\sigma$ and the estimates on the larger domain follow directly.

\begin{figure}[t]
\floatbox[{\capbeside\thisfloatsetup{capbesideposition={right,top},capbesidewidth=5cm}}]{figure}[\FBwidth]
{\caption{An idea of the domain setup for the single fold problem, with the original triangles represented with a dashed outline.}\label{figure3}}
{\includegraphics[width=4.5cm]{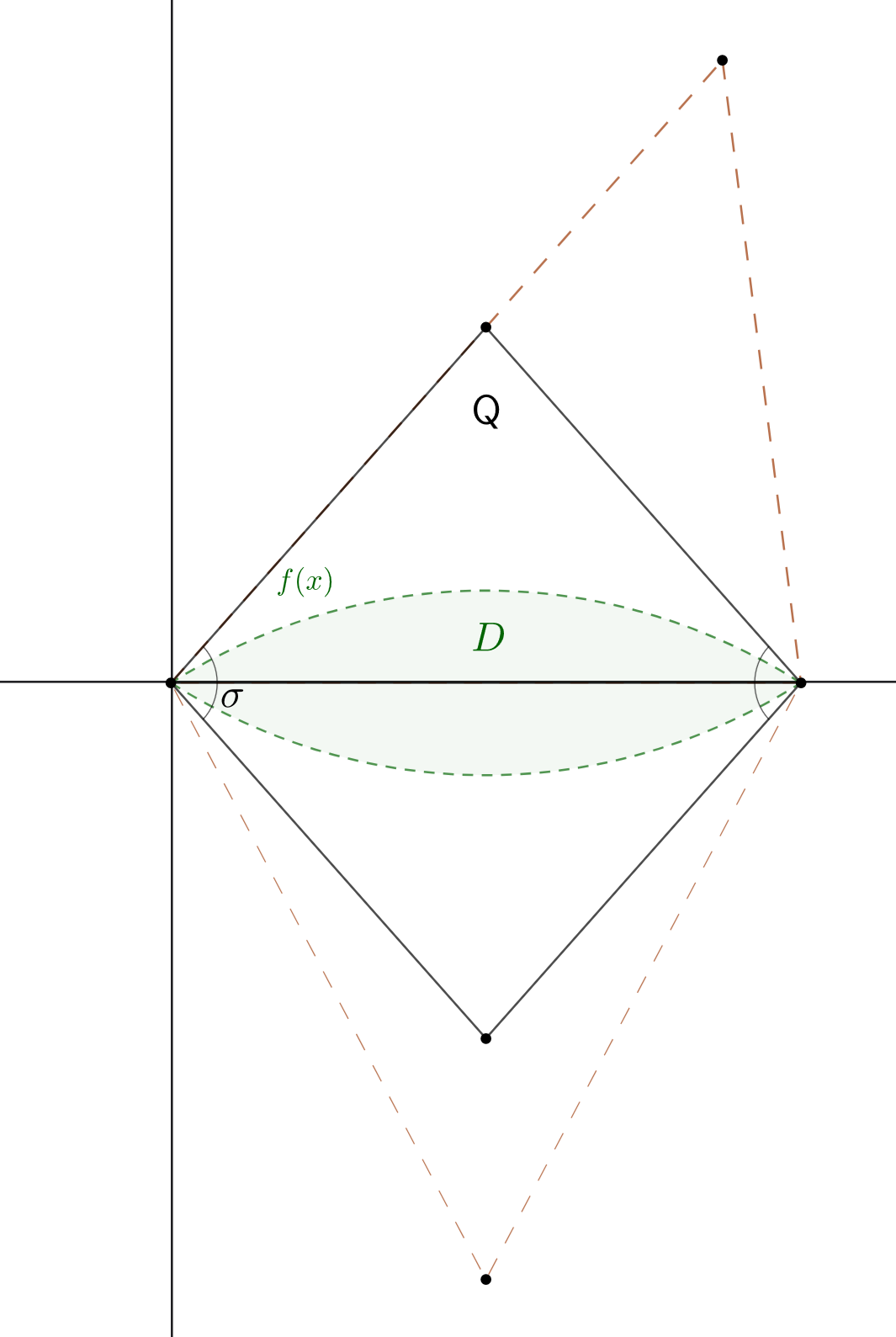}}
\end{figure}


\subsection{Structure of the modified displacement $(w_\sigma,v_\sigma)$}

Let us take a profile function $f\in \mathcal{C}^\infty([0,l];\mathbb{R})$ with $0\le f(x_1)\le \tau \min\{x_1,l-x_1\}$ of a region $D\coloneqq\{x\in Q : \abs{x_2}<f(x_1)\}$ around our fold. We intend to modify our piecewise affine map $(w,v)$ in this region to get a smooth displacement. Let us define:

\begin{equation}\label{vsigma}
    \tilde{v}(x)=a\cdot x + \lambda f(x_1)\gamma\left(\frac{x_2}{f(x_1)}\right)
\end{equation}
\begin{equation}\label{wsigma}
    \tilde{w}(x)=
    Bx-\lambda f(x_1) \gamma\left(\frac{x_2}{f(x_1)}\right) a-\frac{1}{2}\lambda^2 f(x_1)\psi\left(\frac{x_2}{f(x_1)}\right)  e_2+\beta(x)e_1
\end{equation}

where we take $\gamma,\psi \in \mathcal{C}^2([-1,1];\mathbb{R})$ to be our connecting functions, so they must satisfy:
\begin{equation} \label{gammabase}
    \gamma(t)=\begin{cases} 0 \quad \textrm{in a neighborhood of } t=1\\
    t \quad \textrm{in a neighborhood of } t=-1
    \end{cases}
\end{equation}
\begin{equation}\label{chibase}
    \psi (t)=\begin{cases} 0 \quad \textrm{in a neighborhood of } t=1\\
    t \quad \textrm{in a neighborhood of } t=-1
    \end{cases}
\end{equation}

And we take a correction $\beta\in \mathcal{C}^2(\mathbb{R}^2;\mathbb{R})$ in such a way that 
\begin{equation}\label{betabase}
    \beta(x)=0\quad \textrm{ for }\quad \abs{x_2}\simeq f(x_1) 
\end{equation}

We can see that the newly constructed function coincides with $(w,v)$ on the border of $D$. Our displacement function will be then be defined as such:
\begin{equation}
    (w_\sigma, v_\sigma)(x)=\begin{cases}
        (w,v)(x) &\textrm{for }x\in Q\setminus D\\
        (\tilde{w},\tilde{v})(x) &\textrm{for }x\in D
    \end{cases}
\end{equation}
If we can construct the above mentioned functions we may build a displacement function that only differs from our original piecewise affine map $(w,v)$ in the region $D$; also, for $x$ in a neighborhood of $\partial D$ we would have $x_2\simeq \pm f(x_1)$ and by \eqref{gammabase}, \eqref{chibase} and \eqref{betabase} we have then $(\tilde{w},\tilde{v})(x)=(w,v)(x)$. Since all functions taken in the definition of $(w_\sigma, v_\sigma)(x)$ are of class at least $\mathcal{C}^2$, continuity of the definition in a neighborhood of the border of $D$ is enough to have a deformation that is indeed $\mathcal{C}^2$ away from the vertices as required. In order to prove Theorem \ref{mainthm2} it then suffices to prove that there exist some functions $f,\gamma,\psi,\beta$ as above and such that $(w_\sigma,v_\sigma)$ satisfies the estimates from Lemma \ref{singlefoldthm}.

We may start by estimating the stretching component of the energy for our fold: $\frac{1}{2} \nabla v_\sigma\otimes \nabla v_\sigma + \sym \nabla w_\sigma$. We already know this value to be 0 outside of $D$ so we will restrict ourselves to $D$. We define $p(t)=\gamma(t)-t\gamma'(t)$ and $q(t)=\psi(t)-t\psi'(t)$ for all $t\in[-1,1]$. Firstly we have:

\begin{gather}\label{vtilde}
    \nabla \tilde{v}=a+\lambda f'pe_1 + \lambda\gamma' e_2\\
    \label{wtilde} \nabla \tilde{w}=B-\lambda f' p(a\otimes e_1)-\lambda\gamma'(a\otimes e_2)-\frac{\lambda^2}{2}(f' q e_2\otimes e_1 +\psi' e_2\otimes e_2)+\beta_1  e_1\otimes  e_1 + \beta_2  e_1\otimes  e_2
\end{gather}
where to simplify our notation we have set $\beta_{i}\coloneqq\partial_{x_{i}}\beta(x)$ for $i\in \{1,2\}$ and we assume $f,f'$ to always be computed in $x_1$ (unless otherwise specified), and $\gamma,\psi,p,q$ to similarly always be computed in the variable $t=x_2/f(x_1)$. Hence we get:

\begin{equation}
\begin{split}\label{derivv}
    \frac{1}{2} \nabla \tilde{v}\otimes\nabla\tilde{v}=&\frac{1}{2}[a\otimes a + \lambda f'p( e_1\otimes a+a\otimes e_1) +\lambda\gamma'(a\otimes  e_2+ e_2\otimes a)+\\&+\lambda^2f'^2p^2 e_1\otimes e_1 + \lambda^2 f'\gamma'p( e_1\otimes e_2+ e_2\otimes e_1)+\lambda^2\gamma'^2 e_2\otimes e_2]
\end{split}
\end{equation}
\begin{equation}\label{derivw}
    \begin{split}
    \sym \nabla\tilde{w}=&\sym B-\frac{\lambda}{2}f'p( e_1\otimes a +a\otimes e_1)-\frac{\lambda}{2}\gamma'(a\otimes e_2+ e_2\otimes a)+\\&+\beta_1 e_1\otimes e_1+\left(\frac{\beta_2}{2}-\frac{\lambda^2}{4}f'q\right)( e_1\otimes e_2+ e_2\otimes e_1)-\frac{\lambda^2}{2} \psi' e_2\otimes e_2
\end{split}
\end{equation}
By combining \eqref{vonkarman2} with \eqref{derivv} and \eqref{derivw} we finally get:
\begin{equation}\label{stretching1}
\begin{split}
    \frac{1}{2} \nabla \tilde{v}\otimes\nabla\tilde{v}+\sym \nabla\tilde{w}= \left(\frac{\lambda^2}{2}f'^2p^2+\beta_1\right) e_1\otimes e_1 + \left(\frac{\lambda^2}{2} f'\zeta+\frac{\beta_2}{2}\right)( e_1\otimes e_2+ e_2\otimes e_1)+\\+\frac{\lambda^2}{2}(\gamma'^2-\psi') e_2\otimes e_2
    \end{split}
\end{equation}
where we have defined $\zeta(t)=\gamma'(t)p(t)-\frac{1}{2}q(t)$.


\subsection{Choice of $\gamma,\psi,\beta$ and $f$.}

While we do not want the stretching energy to be zero, so that the interplay between stretching and bending energy may provide an optimal upper bound, we will try to simplify \eqref{stretching1} by choosing the proper functions $\gamma,\psi,\beta$. We will, first of all, construct $\gamma$ and $\psi$ in such a way that $\psi'-\gamma'^2=0$ so that \eqref{stretching1} becomes:
\begin{equation}\label{stretching2}
    \frac{1}{2} \nabla \tilde{v}\otimes\nabla\tilde{v}+\sym \nabla\tilde{w}= \left(\frac{\lambda^2}{2}f'^2p^2+\beta_1\right) e_1\otimes e_1 + \left(\frac{\lambda^2}{2} f'\zeta+\frac{\beta_2}{2}\right)( e_1\otimes e_2+ e_2\otimes e_1)
\end{equation}
To do so, it suffices to define $\psi(t)=\int_{-1}^{t}\gamma'^2(s)ds-1$, where if $\gamma\in \mathcal{C}^2([-1,1];\mathbb{R})$ satisfies \eqref{gammabase} and $\int_{-1}^1\gamma'^2=1$, then \eqref{chibase} holds as well.

It remains to define $\gamma$ so that it satisfies the properties mentioned above. We will proceed by taking:
\[\gamma_0(t)=\begin{cases}
    0 \quad&\textrm{for } t\in[0,1]\\
    t \quad&\textrm{for } t\in[-1,0]
\end{cases}\]
then let $\rho_{1/k}\in \mathcal{C}^\infty$ be an even bump function such that $\int\abs{\rho_{1/k}}=1, \rho_{1/k}\ge 0$ and $\supp (\rho_{1/k})\subset\left[-\frac{1}{k},\frac{1}{k}\right]$.

Then (by extending $\gamma_0(t)=t$ for $t<0$) $\gamma_1(t)=\gamma_0*\rho_{1/3}(t)$ is such that $\gamma_1'=\gamma_0'*\rho_{1/3}(t)$, and since $\gamma'_1\in [0,1]$, we get:
\begin{equation}\label{norm}
    \int_{-1}^{1}|\gamma'_1|^2\le \int_{-1}^1 \gamma_1' = \gamma_1(1)-\gamma_1(-1)= 1
\end{equation}
Then it suffices to define, for some $\epsilon\ge0$ to be specified later on:
\begin{equation}\label{gammafinal}
    \gamma(t)\coloneqq\gamma_1(t)+\epsilon\rho_{1/6}(t-1/2)
\end{equation}
hence:
\[
    \gamma(t)=\begin{cases}
        \gamma_1(t) & \textrm{for } t\in\left[-1,\frac{1}{3}\right]\cup \left[\frac{2}{3},1\right]\\
        \epsilon\rho_{1/6}(t-1/2) & \textrm{for } t\in\left[\frac{1}{3},\frac{2}{3}\right]
    \end{cases}
\]
We may observe that \eqref{gammabase} holds and also, by construction, $\gamma$ is smooth. Lastly we may see that for $\epsilon=0$ we have by \eqref{norm}: $\int_{-1}^1\abs{\gamma'}^2=\int_{-1}^1\abs{\gamma'_1}^2\le 1$. On the other hand:
\[
    \norm{\gamma'}_{L^2}^2 = \int_{-1}^1\left[\abs{\gamma'_1(t)}^2+2\epsilon\gamma'_1(t)\rho'_{1/6}(t-1/2)+\epsilon^2\rho'^2_{1/6}(t-1/2)\right]dt=\norm{\gamma_1'}_{L^2}^2+\epsilon^2\lVert\rho'_{1/6}\rVert^2_{L^2}
\]
Since $\lVert{\rho'_{1/6}}\rVert_{L^2}>0$, as $\epsilon\to+\infty$ we get $\norm{\gamma'}_{L^2}^2\to+\infty$ then by continuity there must exist some $\epsilon>0$ such that $\norm{\gamma'}^2_{L^2}=1$; we may then choose $\epsilon$ as such, thus \eqref{chibase} is satisfied.\\

We now turn to choosing $\beta$. We would like to further simplify \eqref{stretching1} by having $\lambda^2f'\zeta+\beta_2=0$. We can choose $\beta(x)=-\lambda^2f(x_1)f'(x_1)\omega\left(\frac{x_2}{f(x_1)}\right)$, with $\omega(t)\coloneqq \int_{-1}^t\zeta(s) ds$. We have:
\[
    \beta_2 =-\lambda^2 f' \zeta
\]
as we wanted. We may also see that $\beta$ is $\mathcal{C}^2$ as long as $\zeta$ is at least $\mathcal{C}^1$ (which is satisfied by our choice of $\gamma$ and $\psi$). We may also see that since $p,q=0$ around $\pm 1$ then $\zeta(t)=0$ for $t$ in a neighborhood of $\pm 1$, also via a triangular integration:
\[
    \int_{-1}^1 q(t)dt=\int_{-1}^1\left(\int_{-1}^t\gamma'^2(s)ds-1-t\gamma'^2(t)\right)dt=\int_{-1}^1((1-2t)\gamma'^2(t)-1) dt=-1-2\int_{-1}^1 t\gamma'^2(t)dt
\]
And consequently:
\[\begin{split}
    \int_{-1}^1 \zeta(t)dt=&\int_{-1}^1\gamma'(t)p(t)dt-\frac{1}{2}\int_{-1}^1q(t)dt=\\
    =&\int_{-1}^1\gamma'(t)(\gamma(t)-t\gamma'(t))dt+\frac{1}{2}+\int_{-1}^1 t\gamma'^2(t)=\\
    =&\frac{1}{2}\left( \gamma(1)^2-\gamma(-1)^2\right)+\frac{1}{2}=0
    \end{split}
\]
As such $\omega(t)=0$ in a neighborhood of $\pm 1$ and \eqref{betabase} is satisfied. A straightforward computation gives us:
\[
    \beta_1(x)=\lambda^2\left[f'^2(x_1)\xi\left(\frac{x_2}{f(x_1)}\right)-f(x_1)f''(x_1)\omega \left(\frac{x_2}{f(x_1)}\right)\right]
\]
where we have defined $\xi(t)\coloneqq t\zeta(t)-\omega(t)$. Now we turn to choose $f$. We will make direct use of Lemma 2.4 from \cite{contimaggi}.

\begin{lemma}
    For any choice of $\tau, l >0$, $\sigma \in (0,l/8)$ there exists a function $f\in \mathcal{C}^\infty([0,l])$ such that:
    \begin{enumerate}
        \item[(i)] $f$ is even with respect to $x=l/2$, or in other words $f(x)=f(l-x)$;
        \item[(ii)] $f$ and its first three derivatives are controlled on $(0,l/2)$ by $f_0(x)=\tau \sigma^{1/3}(x+\sigma)^{2/3}-\tau\sigma$, in the sense that there is some $C>0$ such that for all $x\in(0,l/2)$:
        \begin{equation}\label{fbound1}
            \frac{1}{C}f_0\le f \le f_0, \quad \abs{f'}\le C\abs{f_0'}, \quad \abs{f''}\le C\abs{f_0''}, \quad \abs{f'''}\le C\abs{f_0'''}
        \end{equation}
    \end{enumerate}
\end{lemma}
We also note that:
\begin{equation}\label{fbound2}
    f'_0(x)=\frac{2}{3} \tau\sigma^{1/3}(x+\sigma)^{-1/3},\quad
    f''_0(x)=-\frac{2}{9} \tau\sigma^{1/3}(x+\sigma)^{-4/3}
\end{equation}


\subsection{Upper bound scaling on a single fold}

Let us now compute the bound for the stretching energy of the single fold as requested. We have from \eqref{stretching1} after substituting in $\beta$ and letting $t=x_2/f(x_1)$:
\begin{equation*}
\begin{split}
    \iint_{Q_\sigma}&\abs{\frac{1}{2} \nabla \tilde{v}\otimes\nabla\tilde{v}+\sym \nabla\tilde{w}}^2dx =
    \\&\quad\quad\quad\quad\quad\quad\quad\quad =\lambda^4\int_{\sigma'}^{l-\sigma'}\int_{-1}^1\left[f'^2(x_1)\left(\frac{p^2(t)}{2}+\xi(t)\right)-f(x_1)f''(x_1)\omega(t) \right]^2 f(x_1) dt dx_1\\
    &\quad\quad\quad\quad\quad\quad\quad\quad\le C_\gamma \lambda^4 \left[\int_\sigma^l f(x_1)f'^4(x_1)dx_1 + \int_\sigma^l f^3(x_1) f''^2(x_1)dx_1\right]
\end{split}
\end{equation*}
with $C_\gamma$ some constant depending only on the choice of $\gamma$, and we have ignored the double product of the square by choosing at most a larger $C_\gamma$. By using \eqref{fbound1} and \eqref{fbound2} we finally get:
\begin{equation}\label{stretchingfinal}
    \iint_{Q_\sigma}\abs{\frac{1}{2} \nabla \tilde{v}\otimes\nabla\tilde{v}+\sym \nabla\tilde{w}}^2dx \le C\lambda^4\tau^5\int_{\sigma}^{\frac{l}{2}}\frac{\sigma^{5/3}}{(x_1+\sigma)^{2/3}}dx_1\le C\lambda^4\tau^5l^{1/3} \sigma^{5/3}
\end{equation}
Where, after integrating, we have ignored higher powers of $\sigma$, since we may assume that we are computing this bound for $\sigma$ sufficiently small. 

We may now turn to bound the bending energy $h^2\int_{Q_\sigma} \abs{\nabla^2 v_\sigma}^2$. We remark that $\abs{x_2}\le f(x_1)$ on $Q_\sigma$ also $p'(t)=-t\gamma''(t)$. Then, once again using $t=x_2/f(x_1)$:
\[
    \nabla^2v_\sigma (x)=\lambda \left[\left(f''p+\frac{f'^2}{f} t^2\gamma''\right) e_1\otimes e_1-\frac{f'}{f} t\gamma''( e_1\otimes e_2+ e_2\otimes e_1)+\frac{\gamma''}{f} e_2\otimes e_2\right]
\]
Since $t, p,\gamma', \gamma''$ are uniformly bounded we get, through straightforward algebraic bounds:
\begin{equation}\label{bendingfinal}
    \begin{split}
        \iint_{Q_\sigma}\abs{\nabla^2v_\sigma(x)}^2 dx_1dx_2 &\le C \int_\sigma^l \int_{-1}^1 \lambda^2f \left[\left(f''p+\frac{f'^2}{f} t^2\gamma''\right)^2+\frac{f'^2}{f^2} t^2\gamma''^2+\frac{\gamma''^2}{f^2}\right]dx_1 dt\\
        &\le C\lambda^2 \int_\sigma^l \left(f(x_1)f''^2(x_1)+\frac{f'^4(x_1)+1}{f(x_1)}\right)dx_1\\
        &\le C\lambda^2\int_\sigma^\frac{l}{2} \left(\frac{\tau^3\sigma}{(x_1+\sigma)^2}+\frac{\tau^3\sigma^{5/3}}{(x_1+\sigma)^{8/3}}+\frac{\sigma^{-1/3}}{\tau(x_1+\sigma)^{2/3}}\right)dx_1\\
        &\le C\lambda^2\left(\tau^3+\frac{l^{1/3}}{\tau\sigma^{1/3}}\right)
    \end{split}
\end{equation}
So by merging \eqref{stretchingfinal} and \eqref{bendingfinal} and seeing that we might as well impose $\tau\le 1$:
\begin{equation}\label{sigmafinal}
    I^h[(w_\sigma, v_\sigma), Q_\sigma]\le C\left[\abs{\lambda}^4l^{1/3} \sigma^{5/3}+\abs{\lambda}^2h^2\frac{l^{1/3}}{\tau\sigma^{1/3}}\right]
\end{equation}

We can also see that the desired estimates on the derivatives hold on $Q_h$ by performing direct computations on \eqref{vtilde} and \eqref{wtilde}:
\begin{equation}\label{vbound2}
    \abs{\partial_{x_1}v_\sigma-\partial_{x_1}v}= \abs{\lambda f'p}\le C\tau\abs{\lambda}
\end{equation}
\begin{equation}\label{vbound3}
    \abs{\partial_{x_2}v_\sigma-\partial_{x_2}v}=\abs{\lambda\gamma'}\le C\abs{\lambda}
\end{equation}
\begin{equation}\label{wbound2}
    \abs{\partial_{x_1}w_\sigma-\partial_{x_1}w}\le C\abs{\lambda f'pa+\beta_1 e_1+\frac{1}{2} \lambda^2 f'q e_2}\le C\tau(\abs{\lambda}\abs{a} + \lambda^2)
\end{equation}
\begin{equation}\label{wbound3}
    \abs{\partial_{x_2}w_\sigma-\partial_{x_2}w}\le \abs{\lambda \gamma'a-\beta_2 e_1+\frac{1}{2}\lambda^2\psi' e_2}\le C(\abs{\lambda}\abs{a}+\lambda^2)
\end{equation}
\begin{equation}\label{hessianbound}
    \norm{\nabla^2 v_\sigma}_{L^\infty}\le C\frac{\abs{\lambda}}{\tau \sigma}
\end{equation}

and consequently from \eqref{vbound3} and \eqref{wbound3}, on any $x\in Q_h$ it holds:
\begin{equation}\label{vbound1}
    \abs{v_\sigma(x)-v(x)}\le \abs{x_2-f(x_1)}\sup_{Q_h}\abs{\partial_{x_2} v_\sigma-\partial_{x_2} v}\le C\abs{\lambda}\tau \sigma^{1/3} l^{2/3}
\end{equation}
\begin{equation}\label{wbound1}
    \abs{w_\sigma(x)-w(x)}\le \abs{x_2-f(x_1)}\sup_{Q_h}\abs{\partial_{x_2} w_\sigma-\partial_{x_2} w}\le C(\abs{\lambda}+\abs{a})\abs{\lambda}\tau \sigma^{1/3} l^{2/3}
\end{equation} 

which is sufficient to prove Lemma \ref{singlefoldthm}.\qed

\subsection{Upper bound scaling on the full domain}
We will now extend this upper bound to the whole domain and prove Theorem \ref{mainthm2}. Let $(\mathcal{T},(w,v))$ be our piecewise affine map on $\Omega$ with respect to the triangulation $\mathcal{T}$ where we have $\mathcal T=\{T:T\cap\Omega\neq\emptyset\}$ with respect to $\tilde{\mathcal{T}}_l$ from Section 2; let $\nu$ be the maximum number of triangles intersecting at any vertex in $\mathcal{T}$, and let $l_{\max}, l_{\min}$ be respectively the maximum and minimum length among the sides of all triangles and $\tau_{\min}$ the minimum slope. Let $\lambda_*$ be the norm of maximum jump of $\nabla v$, which we may assume to be nonzero (otherwise the function is fully affine and the estimates are trivially true). We notice that there exists some $\sigma^*$ such that for all $x\in \Omega$, $B(x,2\sigma^*)$ intersects at most $2\nu$ triangles of $\mathcal{T}$. Let:
\begin{equation}\label{hcond}
    \sigma<\min\{\sigma^*,l_{\min}/8\}
\end{equation}

We now divide each triangle in $\mathcal{T}$ in three different triangles by connecting the barycenter to the vertices, and then apply Lemma \ref{singlefoldthm} on each pair of adjacent "subtriangles" whose side in common belonged to the original triangulation $\mathcal{T}$. As we now replace $(w,v)$ on each triangles with $(w_\sigma,v_\sigma)$ as we have seen, by construction the definitions of our new deformation are compatible on the whole domain, save for a neighborhood of the vertices of $\mathcal{T}$. Hence, if we call $\Omega_\sigma\coloneqq\bigcup\limits_{i=1}^N B(z_i,\sigma) $ with $z_i$ vertices of $\mathcal{T}$ for $i=1,\dots, N$, then we get from \eqref{sigmafinal} on $\Omega\setminus\Omega_\sigma$:
\begin{equation}\label{finalbound1}
    I^h[(w_\sigma,v_\sigma),\Omega\setminus\Omega_\sigma]\le C\abs{\mathcal{T}}\left[\lambda_*^4l_{\max{}}^{1/3} \sigma^{5/3}+\lambda_*^2h^2\frac{l_{\max{}}^{1/3}}{\tau_{\min{}}\sigma^{1/3}}\right]
\end{equation}
where we have left $(w_\sigma,v_\sigma)=(w,v)$ on the outer edges.


We will then show that the energy estimate on $\Omega_\sigma$ can be made negligible. Let us consider the values of the piecewise affine $\nabla v$ on $\Omega_{2\sigma}$: since it takes at most $2\nu$ different values around each vertex $z_i$ we may choose one of them and call it $F_i$ for $i=1,\dots, N$. We may call $G_i$ the value of $\nabla w$ on the same triangle on which $\nabla v=F_i$ for each vertex $z_i$. We have, by applying Lemma \ref{singlefoldthm} (and more specifically \eqref{2.7} and \eqref{2.8}) on each triangle, and looking at the bounds used in \eqref{bendingfinal}:
\begin{equation}\label{vwbound2}
    \lVert \nabla v_\sigma-\nabla v\rVert_{L^\infty(\Omega\setminus \Omega_\sigma)} \le C_{\tau_{\min}}\lambda_*
\end{equation}
\begin{equation}\label{hessianboundfinal}
	h\lVert \nabla^2v_\sigma \rVert_{L^\infty(\Omega\setminus \Omega_\sigma)}\le C\frac{\lambda_*h}{\sigma}
\end{equation}
\begin{equation}\label{vwbound4}
    \lVert \nabla w_\sigma-\nabla w\rVert_{L^\infty(\Omega\setminus \Omega_\sigma)} \le C_{\tau_{\min}}\lambda_*\norm{\nabla v}_{L^\infty}
\end{equation}
On the other hand, by using \eqref{2.5} and \eqref{2.6} we get:
\begin{equation}\label{vwbound1}
    \lVert v_\sigma-v\rVert_{L^\infty(\Omega\setminus \Omega_\sigma)}\le C_{\tau_{\min}} \lambda_*\sigma^{1/3} l^{2/3}_{\max}
\end{equation}
\begin{equation}\label{vwbound3}
    \lVert w_\sigma-w\rVert_{L^\infty(\Omega\setminus \Omega_\sigma)}\le C_{\tau_{\min}} \lambda_*\norm{\nabla v}_{L^\infty}\sigma^{1/3} l^{2/3}_{\max}
\end{equation}
By applying \eqref{vwbound2} and \eqref{vwbound4} on each $B(z_i,2\sigma)$:
\begin{equation}\label{F}
    \abs{\nabla v_\sigma- F_i}\le \abs{\nabla v_\sigma- \nabla v}+\abs{\nabla v- F_i}\le C_{\tau_{\min}}\nu \lambda_*
\end{equation}
\begin{equation}\label{G}
    \abs{\nabla w_\sigma- G_i}\le \abs{\nabla w_\sigma- \nabla w}+\abs{\nabla w- G_i}\le C_{\tau_{\min}}\nu \lambda_*\norm{\nabla v}_{L^\infty}
\end{equation}
We may now define $(\bar{w},\bar{v}) :\Omega_{2\sigma}\to \mathbb{R}^2\times \mathbb{R}$ defined for all $x\in B(z_i,2\sigma)$ as:
\begin{gather*}
    \bar{v}(x)\coloneqq v(z_i) + F_i \cdot (x-z_i)\\
    \bar{w}(x)\coloneqq w(z_i)+ G_i(x-z_i)
\end{gather*}
It follows that:
\[
    \frac{1}{2} \nabla \bar{v}\otimes\nabla\bar{v}+\sym \nabla \bar{w}=0
\]
and from \eqref{F} and \eqref{G}:
\begin{equation}\label{overline2}
    \lVert \nabla \bar v  - \nabla v_\sigma \rVert_{L^\infty(\Omega_{2\sigma})}\le C \nu \lambda_*, \quad \lVert \nabla \bar w  - \nabla w_\sigma \rVert_{L^\infty(\Omega_{2\sigma})}\le C \nu \lambda_*\norm{\nabla v}_{L^\infty}
\end{equation}
\begin{equation}\label{overline1}
    \lVert  \bar v  -  v_\sigma \rVert_{L^\infty(\Omega_{2\sigma})}\le C \nu \lambda_* \sigma, \quad \lVert  \bar w  -  w_\sigma \rVert_{L^\infty(\Omega_{2\sigma})}\le C \nu \lambda_*\norm{\nabla v}_{L^\infty}\sigma
\end{equation}
Let $B^2(z,r)$ be the ball of center $z$ and radius $r$ in $\mathbb{R}^2$ and let us choose a bump function $\alpha\in \mathcal{C}^\infty_c (B^2(0,2),[0,1])$, with $\alpha(x)=1$ if $x\in B^2(0,1)$. We may build the mollifier:
\[
    \alpha_\sigma(x)\coloneqq \sum_{i=1}^N \alpha\left(\frac{x-z_i}{\sigma}\right)
\]
where we have: 
\[
    \alpha\left(\frac{x-z_i}{\sigma}\right)=\begin{cases}
        1 & \textrm{if } x\in B(z_i,\sigma)\\
        0 & \textrm{if } x\notin B(z_i, 2\sigma)
    \end{cases}
\]
and since $l_{\min}>2\sigma$ the sum is disjoint hence $ \norm{\alpha_\sigma}_\infty=\norm{\alpha}_\infty$ constant. We finally interpolate and define:
\[
    \bar v_\sigma(x)\coloneqq \alpha_\sigma(x) \bar v(x) +(1-\alpha_\sigma(x))v_\sigma(x)
\]
\[
    \bar w_\sigma(x)\coloneqq \alpha_\sigma(x) \bar w(x) +(1-\alpha_\sigma(x))w_\sigma(x)
\]
We will now provide bounds for $\bar v_\sigma$ ($\bar w_\sigma$ is analogous). We may remark that $\abs{\bar v_\sigma-v}\le \alpha_\sigma\abs{\bar{v}-v_\sigma}+\abs{v_\sigma-v}$, hence by using \eqref{vwbound1} and \eqref{overline1}:
\begin{equation}\label{vfinalbound1}
    \norm{\bar v_\sigma-v}_{L^\infty(\Omega)}\le C\left(\nu \sigma\lambda_*\norm{\alpha_\sigma}_{L^\infty(\Omega)} + \lambda_*\sigma^{1/3} l^{2/3}_{\max}\right)\le C_{\Gamma(\mathcal{T})}\lambda_* \sigma^{1/3}
\end{equation}
where $C_{\Gamma(\mathcal{T})}$ depends only on $\nu$ and $\tau_{\min}$ and, as such, it only depends on the degeneracy factor $\Gamma(\mathcal{T})$. On the other hand, $\nabla \bar v_\sigma-\nabla v_\sigma=\alpha_\sigma(\nabla \bar v - \nabla v_\sigma)+(\bar v-v_\sigma)\otimes\nabla\alpha_\sigma$, hence by combining \eqref{vwbound2}, \eqref{overline2} and \eqref{overline1}:
\begin{equation}\label{vfinalbound2}
\begin{split}
    \lVert\nabla \bar v_\sigma - \nabla v\rVert_{L^\infty(\Omega)}&\le \lVert\nabla \bar v_\sigma - \nabla v_\sigma\rVert_{L^\infty(\Omega)}+\lVert\nabla  v_\sigma - \nabla v\rVert_{L^\infty(\Omega)}\le\\
    &\le \lVert \alpha\rVert_{L^\infty}\left[\lVert\nabla \bar v - \nabla v_\sigma\rVert_{L^\infty(\Omega)}+\frac{1}{\sigma}\lVert \bar v -  v_\sigma\rVert_{L^\infty(\Omega)}\right]+\lVert\nabla  v_\sigma - \nabla v\rVert_{L^\infty(\Omega)}\le\\
    &\le C_{\Gamma(\mathcal{T})} \lambda_*
\end{split}
\end{equation}
As for the second derivatives, we have that $\bar v_\sigma$ is affine on $\Omega_\sigma$, and on $\Omega_{2\sigma}$ we have $\abs{\nabla^2 \bar v_\sigma}\le \abs{\nabla^2 v_\sigma}+ \frac{C}{\sigma}\abs{\nabla\bar v-\nabla v_\sigma}+\frac{C}{\sigma^2}\abs{\bar v- v_\sigma}$, hence, by using \eqref{hessianboundfinal}, \eqref{overline2} and \eqref{overline1}:
\begin{equation}\label{vfinalbound3}
    h \lVert \nabla^2 \bar v_\sigma \rVert_{L^\infty(\Omega)}\le C_{\Gamma(\mathcal{T})} \frac{\lambda_*h}{\sigma}
\end{equation}
We may proceed the same way to get, by using \eqref{vwbound4},\eqref{vwbound3},\eqref{overline2},\eqref{overline1}:
\begin{equation}\label{wfinalbound1}
    \lVert \bar w_\sigma - w\rVert_{L^\infty(\Omega)}\le C_{\Gamma(\mathcal{T})}\lambda_*\norm{\nabla v}_{L^\infty} \sigma^{1/3}
\end{equation}
\begin{equation}\label{wfinalbound2}
    \lVert \nabla\bar w_\sigma - \nabla w\rVert_{L^\infty(\Omega)}\le C_{\Gamma(\mathcal{T})}\lambda_*\norm{\nabla v}_{L^\infty}
\end{equation}
All that is left is estimating the energy. For the stretching energy
\[
    S^h[(\bar w_\sigma,\bar v_\sigma),\Omega_{2\sigma}]=\int_{\Omega_{2\sigma}} \abs{\frac{1}{2}\nabla \bar v_\sigma\otimes \nabla \bar v_\sigma+\sym \nabla \bar w_\sigma}^2
\]
    we have, by combining \eqref{vfinalbound2} and \eqref{wfinalbound2}:
\[
    \begin{split}
        S^h[(\bar w_\sigma,\bar v_\sigma),\Omega_{2\sigma}]&\le 
        C\int_{\Omega_{2\sigma}} \abs{\nabla \bar v_\sigma\otimes \nabla \bar v_\sigma-\nabla v\otimes \nabla  v}^2+C\int_{\Omega_{2\sigma}} \abs{\sym \nabla \bar w_\sigma-\sym \nabla  w}^2\\
        &\le C\abs{\Omega_{2\sigma}}\left[ \lVert\nabla \bar v_\sigma - \nabla v\rVert_{L^\infty(\Omega)}^2\lVert\abs{\nabla \bar v_\sigma} + \abs{\nabla v}\rVert_{L^\infty(\Omega)}^2+ \lVert \nabla\bar w_\sigma - \nabla w\rVert_{L^\infty(\Omega)}^2 \right]\\
        &\le C_{\Gamma(\mathcal{T})} \abs{\mathcal{T}} \lambda_*^2\norm{\nabla v}_{L^\infty}^2\sigma^2
    \end{split}
\]
where we have used that $\abs{\Omega_{2\sigma}}\sim \abs{\mathcal{T}} \sigma^2$. By adding the bending energy that we may estimate through \eqref{vfinalbound3}: 
\[
    I^h[(\bar w_\sigma,\bar v_\sigma),\Omega_{2\sigma}]\le C_{\Gamma(\mathcal{T})} (\abs{\mathcal{T}} \lambda_*^2\norm{\nabla v}_{L^\infty}^2\sigma^2 + \abs{\mathcal{T}}\lambda_*^2h^2)\le C_{\Gamma(\mathcal{T})} \abs{\mathcal{T}} \lambda_*^2(\norm{\nabla v}_{L^\infty}^2\sigma^2+h^2)
\]
And finally, by combining it with \eqref{finalbound1} and we get:
\begin{equation}\label{finalboundsigma}
    I^h[(\bar w_\sigma,\bar v_\sigma),\Omega]\le C_{\Gamma(\mathcal{T})} \abs{\mathcal{T}}\left[\lambda_*^4\sigma^{5/3}+\lambda_*^2\frac{h^2}{\sigma^{1/3}}+\lambda_*^2\norm{\nabla v}_{L^\infty}^2\sigma^2\right]
\end{equation}

We may then optimize by choosing $\sigma= h/\lambda_*$ which we may have if $h<\lambda_* l_{\min{}}/8$, as long as in Theorem \ref{mainthm2} we choose $h^*\ll l$. Then, for simplicity, with an abuse of notation we will call the resulting mapping $(w_h,v_h)$ and the union of $h/\lambda_*$ radius balls around the vertices will be referred to as $\Omega_h$; so finally we obtain from \eqref{finalboundsigma}:
\begin{equation*}
    I^h[(\bar w_\sigma,\bar v_\sigma),\Omega]\le C_{\Gamma(\mathcal{T})} \abs{\mathcal{T}} \lambda_*^{7/3}h^{5/3}
\end{equation*}
giving the desired scaling of order $h^{5/3}$ (the $h^2$ term has been absorbed for $h$ small enough). Furthermore \eqref{estimate3.2} and \eqref{estimate3.3} directly follow on $\Omega\setminus\Omega_h$ from \eqref{vbound1} and \eqref{wbound1} and on $\Omega_h$ from \eqref{vfinalbound1} and \eqref{wfinalbound1}. \qed

\section{An example of an obstacle problem}

We have provided with Corollary \ref{positional} an example of an obstacle problem that can be solved through the application of Theorem \ref{mainthm}. We will now provide the proof of Corollary \ref{positional} through Theorem \ref{mainthm}.
\begin{proof}
    Since $v_{\sup}-v_{\inf}>r>0$ it follows that there exists $v_0\in \mathcal{C}^1(\bar{\Omega};\mathbb{R})\cap \mathcal{C}^0(\bar\Omega;\mathbb{R})$ such that:
    \begin{gather*}
        v_{\inf}<v_0<v_{\sup}\\
       \inf_{\Omega} (v_{\sup}-v_0)>\frac{r}{4}\\
        \inf_{\Omega} (v_0-v_{\inf})>\frac{r}{4}
    \end{gather*}
    This follows for example from properly mollifying $\bar{v}=\frac{v_{\sup}+v_{\inf}}{2}$ into some function $v_0$ such that $\norm{v_0-\bar{v}}_\infty<r/4$. 

    Now we can show that for any such $v_0$ there exists $w_0$ such that $(w_0, v_0)$ satisfies \eqref{short}. Then the proof follows by applying Theorem \ref{mainthm} with $r/4$. On the other hand:
    \[
        \nabla v_0\otimes \nabla v_0 \le \norm{\nabla v_0}^2\id
    \]
    It then suffices to take some linear horizontal displacement such as $w_0(x)=-\frac{\norm{\nabla v_0}^2}{2} x$ to verify \eqref{short} and prove our thesis.
\end{proof}

\ \\

\appendix

\section{Transformations for solutions (and short maps) with respect to the von Kármán constraint}

In this appendix we will look at some of several transformations that, when applied to $w$ and $v$, leave the von Kármán constraint unchanged. It is straightforward to see that, for example, a translation of the domain $\Omega$ or of the image leaves unchanged a condition of the following kind:
\begin{equation}\label{hypA1}
        \frac{1}{2} \nabla v \otimes \nabla v+\sym \nabla w \le -\delta\id_2
\end{equation}

and will also preserve distances between different functions since they are isometries. Hence, when working to approximate a solution to \eqref{vonkarman} we may assume without loss of generality that a certain point of the solution is 0 while preserving bounds on the distances of our approximation. We will see some other useful transformations that preserve \eqref{hypA1} in the following sections.

\subsection{Rotations of the domain $\Omega$}
\begin{lemma}\label{FR}
    Let $(w,v)\in W^{1,\infty}(\Omega;\mathbb{R}^2\times\mathbb{R})$ be such that \eqref{hypA1} holds for some $\delta\in\mathbb{R}$. For any $R\in O(2)$ isometry we may define for all $y=R^{-1} x\in R^{-1}\Omega$: 
    \[
        \mathcal{D}_R : (w(x),v(x))\mapsto (w_R(y), v_R(y))=\left(R^T w(Ry) ,v(Ry)\right)
    \]
   Then, for any other $(\bar w,\bar v)\in W^{1,\infty}(\Omega;\mathbb{R}^2\times\mathbb{R})$ we have that $(w_R,v_R)\in W^{1,\infty}(R^{-1}\Omega;\mathbb{R}^2\times\mathbb{R})$ satisfies:
    \begin{gather}\label{FR1bis}
        \frac{1}{2} \nabla v_R \otimes \nabla v_R+\sym \nabla w_R \le -\delta\id_2 \\
        \label{FR2bis}
        \norm{(\bar w_R,\bar v_R)-(w_R,v_R)}_{W^{1,\infty}(R^{-1}\Omega;\mathbb{R}^2\times\mathbb{R})}\le \norm{(\bar w,\bar v)-(w,v)}_{W^{1,\infty}(\Omega;\mathbb{R}^2\times\mathbb{R})}\\
        \label{FR3}
        \mathcal{D}_R^{-1}=\mathcal{D}_{R^{-1}}
    \end{gather}
    More specifically if we have an equality for \eqref{hypA1} in the hypothesis, then \eqref{FR1bis} also holds with an equality.

\end{lemma}

\begin{proof}
    The expressions \eqref{FR2bis} and \eqref{FR3} are straightforward to see, we will focus on \eqref{FR1bis}:
\begin{equation}\label{rotation}
    \begin{split}
        \frac{1}{2} (\nabla v_R\otimes\nabla v_R)(y)+\sym\nabla w_R(y) &= R^T\left(\frac{1}{2}\nabla v^T(Ry)\nabla v(Ry) +\nabla w(Ry)+\nabla w^T(Ry)\right)R=\\
        &=R^T\left(\frac{1}{2}(\nabla v\otimes\nabla v) (x)+\sym \nabla w(x)\right) R\le -\delta R^T R\le -\delta \id_2
    \end{split}
\end{equation}
which satisfies our claim. We may also remark from \eqref{rotation} that the equality case holds if we have an equal sign in \eqref{hypA1}.
\end{proof}

\subsection{A class of quasi-affine transformations}
We will now look at another class of transformations which are almost affine on $(w,v)$ and depend on a given vector $\alpha$:
\begin{lemma}\label{talpha}
    Let $(w,v)\in W^{1,\infty}(\Omega;\mathbb{R}^2\times\mathbb{R})$. For any given vector $\alpha\in \mathbb{R}^2$ we may define: 
    \[
        \mathcal{S}_\alpha : (w(x),v(x))\mapsto (w_{(\alpha)}(x), v_{(\alpha)}(x))=\left(w(x)-v(x)\alpha-\frac{1}{2}(\alpha\otimes \alpha)x, v(x)+\alpha \cdot x\right)
    \]
    Then, for any other $(\bar w,\bar v)\in W^{1,\infty}(\Omega;\mathbb{R}^2\times\mathbb{R})$ we have:
    \begin{gather}\label{talpha0}
        \frac{1}{2} \nabla v_{(\alpha)} \otimes \nabla v_{(\alpha)}+\sym \nabla w_{(\alpha)}=\frac{1}{2} \nabla v \otimes \nabla v+\sym \nabla w \\
        \label{talpha2}
        \norm{\bar v_{(\alpha)}-v_{(\alpha)}}_{L^\infty(\Omega;\mathbb{R})}=\norm{\bar v-v}_{L^\infty(\Omega;\mathbb{R})}\\
        \label{talpha3}
        \norm{\bar w_{(\alpha)}-w_{(\alpha)}}_{L^\infty(\Omega;\mathbb{R}^2)}\le \norm{\bar w-w}_{L^\infty(\Omega;\mathbb{R}^2)}+\abs{\alpha}\norm{\bar v-v}_{L^\infty(\Omega;\mathbb{R})}\\
        \label{talphagrad}
        \abs{v_{(\alpha)}(x)- v_{(\alpha)}(y)}\le \abs{v(x)-v(y)}+\abs{\alpha}\abs{x-y}\\
        \label{talphagrad2}
        \abs{w_{(\alpha)}(x)-w_{(\alpha)}(y)}\le \abs{w(x)-w(y)}+\abs{\alpha}\abs{v(x)-v(y)}+\frac{1}{2}\abs{\alpha}^2\abs{x-y}\\
        \label{talpha4}
        \mathcal{S}_{-\alpha}\mathcal{S}_{\alpha}=\id
    \end{gather}
\end{lemma}
\begin{proof}
    We can see that \eqref{talpha2}, \eqref{talpha3} and \eqref{talphagrad} descend directly from the construction. To verify \eqref{talpha0} we will just proceed by direct computation:
    \[
        \frac{1}{2} \nabla v_{(\alpha)} \otimes \nabla v_{(\alpha)}+\sym \nabla w_{(\alpha)}= \frac{1}{2} \nabla v \otimes \nabla v+\sym(\alpha\otimes \nabla v)+\frac{1}{2}\alpha\otimes\alpha+\sym \nabla w - \sym(\alpha\otimes \nabla v)-\frac{1}{2}\alpha\otimes\alpha
    \]
    which also shows us that the equality case is verified, and that more specifically if \eqref{hypA1} holds for $(w,v)$ then we have:
    \begin{equation}\label{talpha1}
        \frac{1}{2} \nabla v_{(\alpha)} \otimes \nabla v_{(\alpha)}+\sym \nabla w_{(\alpha)} \le -\delta\id_2
    \end{equation}
    Finally, let us compute $\mathcal{S}_{-\alpha}(w_{(\alpha)},v_{(\alpha)})$:
    \[
    \begin{split}
        \mathcal{S}_{-\alpha}(w_{(\alpha)}(x),v_{(\alpha)}(x))&=\left(w_{(\alpha)}(x)+v_{(\alpha)}(x)\alpha-\frac{1}{2}(\alpha\otimes \alpha)x, v_{(\alpha)}(x)-\alpha \cdot x\right)=\\
        &=\left(w(x)-v(x)\alpha-\frac{1}{2}(\alpha\otimes \alpha)x+(v(x)+\alpha \cdot x)\alpha-\frac{1}{2}(\alpha\otimes \alpha)x, v(x)\right)=\\
        &=(w(x),v(x))
    \end{split}
    \]
    This concludes the proof. 
\end{proof}
We could also remark that we have defined a group of transformations $\mathcal{S}_\alpha$ with identity $\mathcal{S}_0$ under which the vK constraint is invariant.

\newpage
\nsection{Acknowledgements}

I will start by thanking my PhD advisor, Professor Mohammad Reza Pakzad for suggesting me this topic and for the many ideas and moments of guidance he gave me to help me push through each and every problem that arose while working on this paper; my parents for the constant help and support through the hardships we have faced in these last few years; and finally one of my closest friends, Giacomo Ficarra, for being a good listener.

\section*{Declarations}

\paragraph{Funding.}
This paper has been written thanks to the PhD funding I have received from the University of Toulon.

\paragraph{Conflict of interest.}
The author has no competing interests to declare that are relevant to the content of this article.

\paragraph{Data availability.}
Data sharing is not applicable to this article, as no datasets were generated or analysed during the current study.

\addcontentsline{toc}{section}{References}
\bibliography{mybib}{}
\bibliographystyle{acm}

\end{document}